\documentclass[10pt,a4paper]{article}
\usepackage[latin1]{inputenc}
\usepackage{anysize}
\usepackage{graphicx}
\marginsize{2.5cm}{2.5cm}{2cm}{1cm}
\usepackage{amsmath}
\usepackage{amsfonts}
\usepackage{amssymb}
\newcommand{\bb}[1]{\mathbb{#1}}

\newcommand{\qed}{\begin{flushright} $\Box$ \end{flushright}}
\numberwithin{equation}{section}
\usepackage{url}
\begin{document}
\title{Knight's Tours in Higher Dimensions}
\author{Joshua Erde \thanks{DPMMS, University of Cambridge}\\ \ \\submitted to Electronic Journal of Combinatorics February 9th 2012}
\date{}
\maketitle
\begin{abstract}
In this paper we are concerned with knight's tours on high-dimensional boards.
Our main aim is to show that on the $d$-dimensional board $[n]^d$, with $n$
even, there is always a knight's tour provided that $n$ is
sufficiently large.

In fact, we give an exact classification of the grids
$[n_1] \times \ldots \times [n_d]$ in which there is a knight's tour. This
answers questions of DeMaio, DeMaio and Mathew, and Watkins.
\end{abstract}
\section{Introduction}
A \emph{knight's tour} on an $n \times m$ chessboard is a traversal of the squares of a chessboard using only moves of the knight to visit each square once. A knight's tour is \emph{closed} if the last move of the tour returns the knight to its starting position; otherwise the tour is \emph{open}. Unless otherwise specified we will only consider closed tours in this paper.\\
\\
In graph theoretical terms we can consider an $n \times m$ chessboard as a grid of $n \times m$ points. We associate with this grid a graph, the \emph{knight's graph} $K(n,m)$, where each point is joined to all points a knight's move away. Equivalently $K(n,m)$ is the graph where $V(G) = \{ (i,j)\,:\, 0\leq i \leq n-1\,,\,0\leq j\leq m-1\}$ and $\big((i,j),(k,l)\big) \in E(G) \Leftrightarrow (i - k, j-l) \in \{ (\pm 1 , \pm 2)\,,\,(\pm 2 ,\pm 1)\}$. So an $n \times m$ tour is precisely a Hamiltonian cycle in $K(n,m)$.\\
\\
Beyond this we can define the knight's graph for higher dimensional chessboards. For a board with dimensions $n_1 \times n_2 ... \times n_r$, we define $G=K(n_1, n_2..., n_r)$ in a similar fashion with $$V(G) = \{(i_1,i_2...,i_r)\,:\,0 \leq i_j \leq n_j -1 \text{ for all } j\}$$
\begin{align*}
E(G) = \{\big( (a_1,a_2...,a_r),(b_1,b_2...,b_r) \big)\,:\, \text{ there exists } i_1,i_2 \text{ such that }&\\
|a_{i_1} - b_{i_1}| = 1\,,\,|a_{i_2} - b_{i_2}| = 2 \text{ and } a_i = b_i \text{ for all } i \neq i_1, i_2&\}
\end{align*}\\
\ \\
The question of the existence of knight's tours has been studied by mathematicians through the ages, both professional and amateur. An early solution on the $8 \times 8$ board was found by De Moivre in the 18th century. More recently, Schwenk $\cite{S1991}$ proved\\
\ \\
{\bf Theorem 1 $\cite{S1991}$ :} \\
An $n \times m$ $(n \geq m)$ tour exists if and only if the following conditions hold:\\
1) $n$ or $m$ is even;\\
2) $m \not\in \{1,2,4\}$;\\
3) $(n,m)$ $\neq (4,3)\,,\,(6,3)$ or $(8,3)$.\\
\ \\
Stewart $\cite{S1971}$ constructed some examples of $3$-dimensional knight tours and DeMaio and Mathew $\cite{DM2011}$ fully classified the $3$-dimensional boards which admit knight's tours. They showed\\
\ \\
{\bf Theorem 2 $\cite{DM2011}$ : }\\
A $p \times q \times r$ $(p \geq q \geq r)$ tour exists if and only if the following conditions hold:\\
1) $p,q$ or $r$ is even;\\
2) $p \geq 4$;\\
3) $q \geq 3$.\\ 
\ \\
In the same paper they asked about higher dimensional tours, This question was also asked by DeMaio $\cite{D2007}$ and Watkins $\cite{W2004}$.\\
\ \\
The main result of the paper will be to show\\
\ \\
{\bf Theorem 3 :}\\
For $r \geq 3$ an $n_1 \times n_2 .... \times n_r$ $(n_1 \geq n_2 ... \geq n_r)$ tour exists if and only if the following conditions hold:\\
1) Some $n_i$ is even;\\
2) $n_1 \geq 4$;\\
3) $n_2 \geq 3$.\\
\ \\
We prove this result in Section $2$. In Section $3$ we will consider the problem of knight's tours with more general moves and make several conjectures.\\
\ \\
\ \\
\ \\
\section{Knight's tours in higher dimensions}
In this section we will prove Theorem 3. The proof is inductive on the dimension of the chessboard. However, a slightly stronger hypothesis is needed to complete the induction step which will motivate the definition of a site and a bi-sited tour which follow.\\
\ \\
Given an $n \times m$ tour we call a pair of edges in the tour a \emph{site} if both endpoints of the two edges are two squares away from each other, more precisely, that is two edges $\left((a_1,b_1),(a_2,b_2)\right)$ and $\left((c_1,d_1),(c_2,d_2)\right)$ such that $(|a_1-c_1|,|b_1-d_1|)$ and $(|a_2-c_2|,|b_2-d_2|)$ $\in \{ (0,2),(2,0)\}$. Below are three examples:
\begin{center}
\includegraphics[scale=0.2]{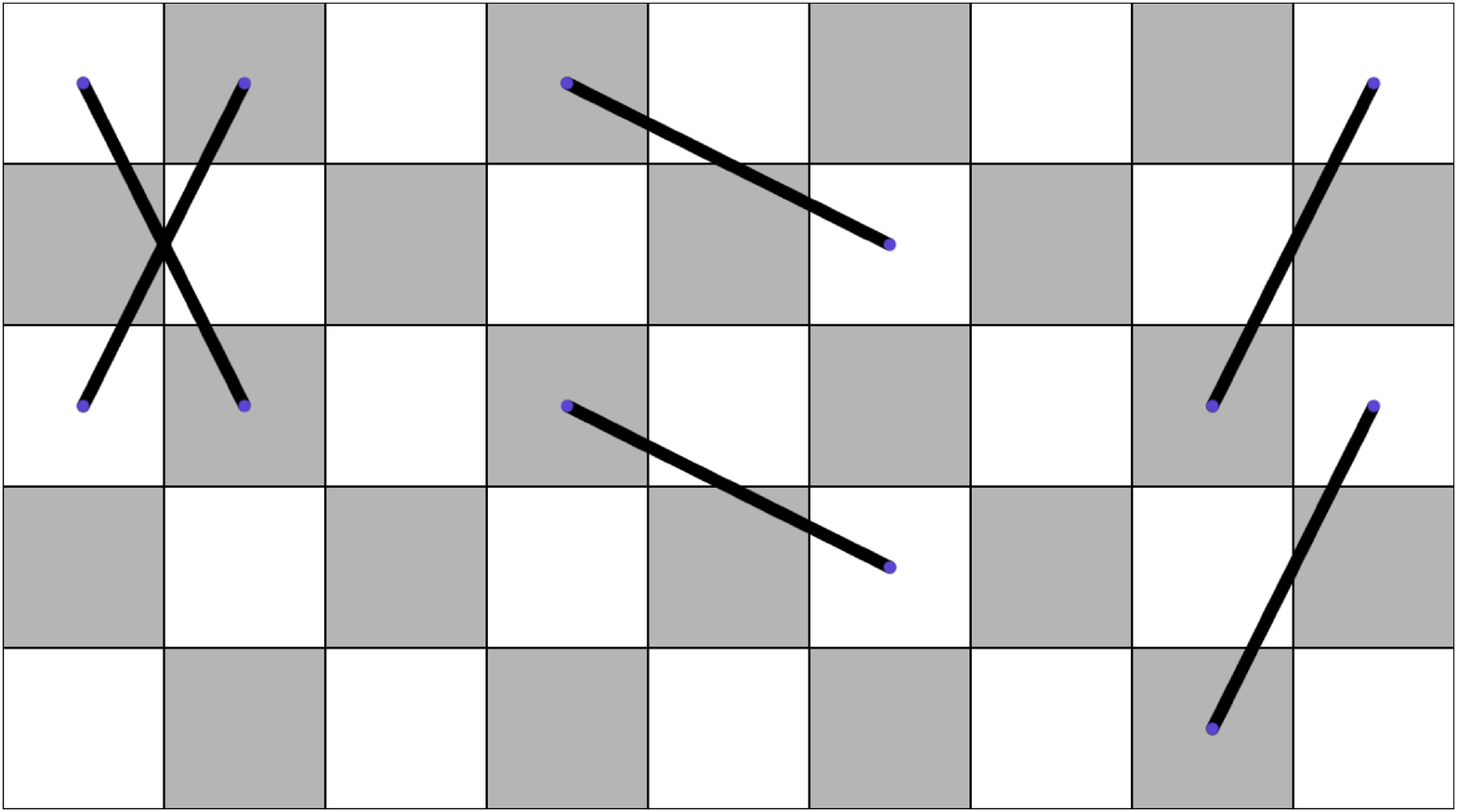}
\end{center}
\ \\
If some $n \times m$ tour contains one of these patterns we can construct an $n \times m \times 2$ tour by placing two copies of the $n \times m$ tour on top of each other, deleting one edge in the site from the top copy and the other edge in the site from the bottom copy and joining up the pairs of vertices which are a knights move apart.\\
\ \\
For example, to construct a $5 \times 6 \times 2$ tour we would place these two tours on top of each other, remove the highlighted edges and add in the edges $\big( (0,2,0)\,,\,(0,4,1) \big)$ and $\big( (1,4,0)\,,\,(1,2,1) \big)$.
\begin{center}
\includegraphics[scale=0.2]{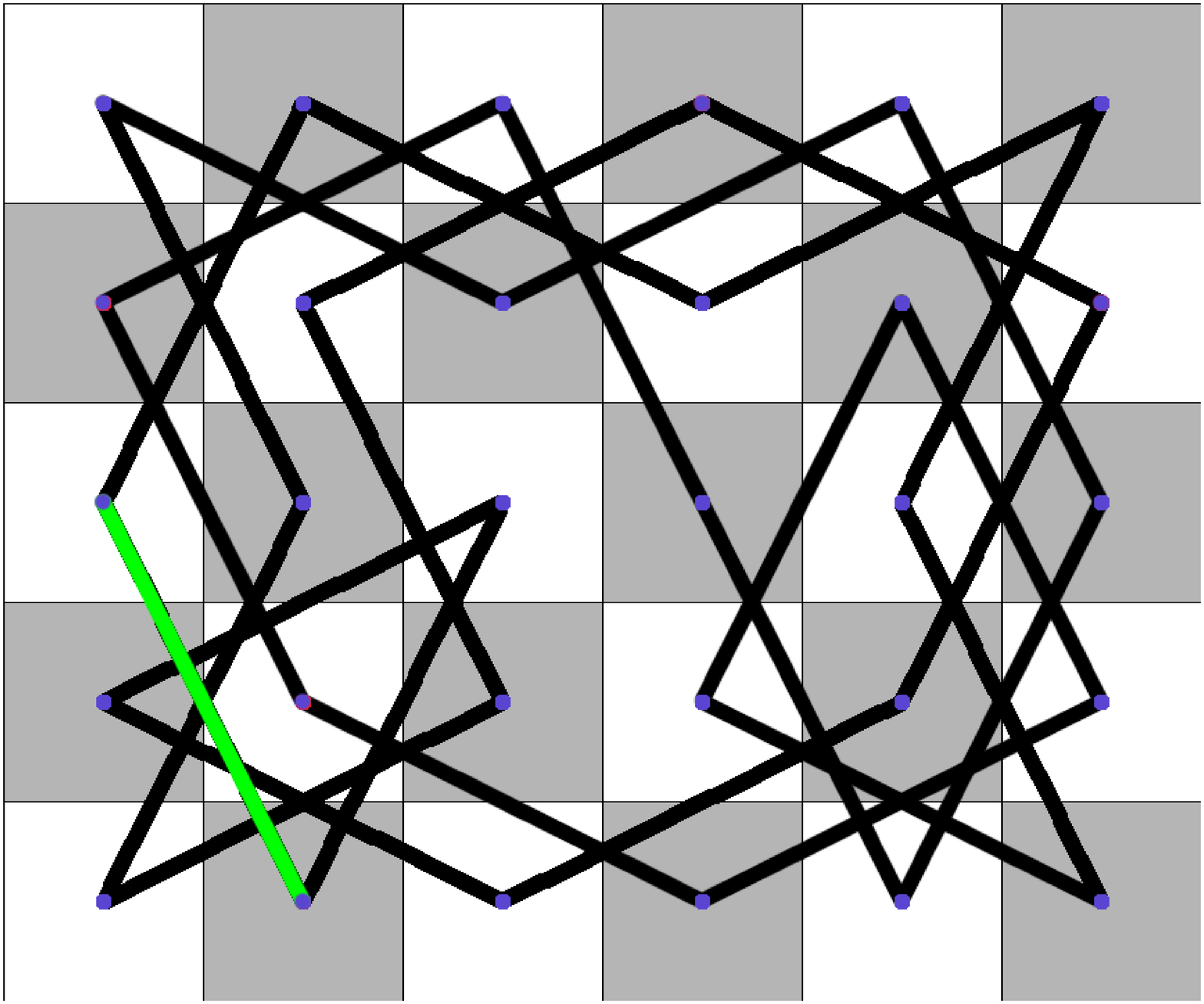}
\includegraphics[scale=0.2]{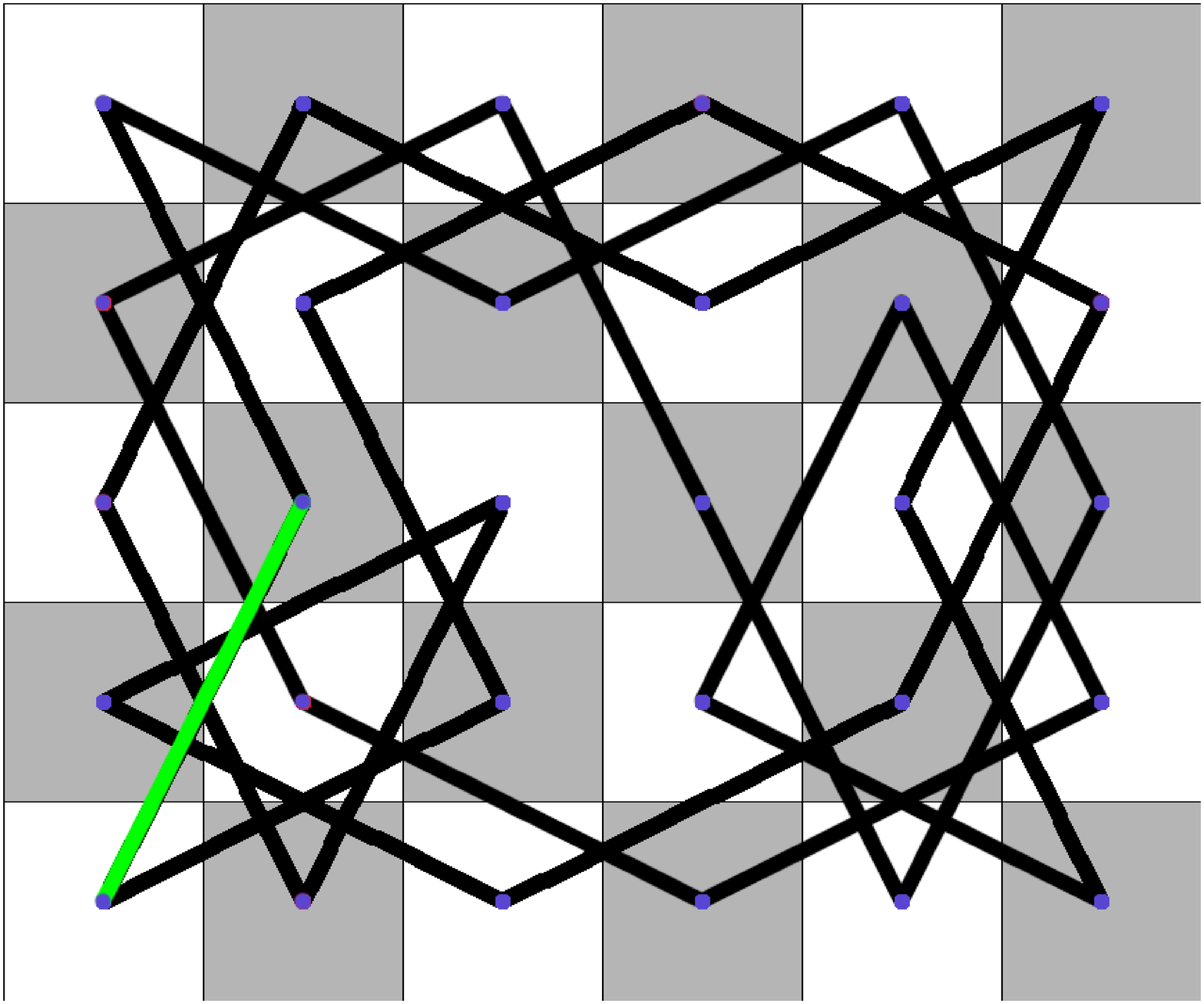}
\end{center}
We call any knight's tour of a chessboard containing two edge disjoint sites \emph{bi-sited}. This is the key idea of the paper, enabling the induction to work.\\
\ \\
{\bf Theorem 4 :}\\
If a bi-sited $n_1 \times .... \times n_r$ tour exists then a bi-sited $n_1 \times .... \times n_r \times p$ tour exists for all $ p \in \bb{N}$.\\
\ \\
{\bf Proof :}\\
We take $p$ copies of the bi-sited $n_1 \times .... \times n_r$ tour and place them on top of each other. We join the first copy to the second copy by the process described above using the first site on both copies, then the second to the third using the second site and so on, alternating sites, until we have formed a $n_1 \times .... \times n_r \times p$ tour. Now there will still be two sites we have not altered during this process, one in the top copy of $n_1 \times .... \times n_r$ and one in the bottom, and so this tour is also bi-sited. \qed
\ \\
{\bf Corollary 5 :}\\
If an $n \times m$ tour exists then so does an $n \times m \times p_1 ....\times p_r$ tour for any $r$ and any $p_1, ..., p_r$.\\
\ \\
{\bf Proof :}\\
Notice that in $K(n,m)$ the vertex $(0,0)$ has degree $2$ and hence any Hamiltonian cycle must contain both edges adjacent to that point, the lines $\big((0,0),(1,2)\big)$ and $\big((0,0),(2,1)\big)$.\\
\ \\
Similarly of the four edges adjacent to the point $(0,2)$, at least $2$ must be in the tour, but 3 of them form sites with the two forced lines, as in the picture below:\\
\ \\
\begin{center}
\includegraphics[scale=0.2]{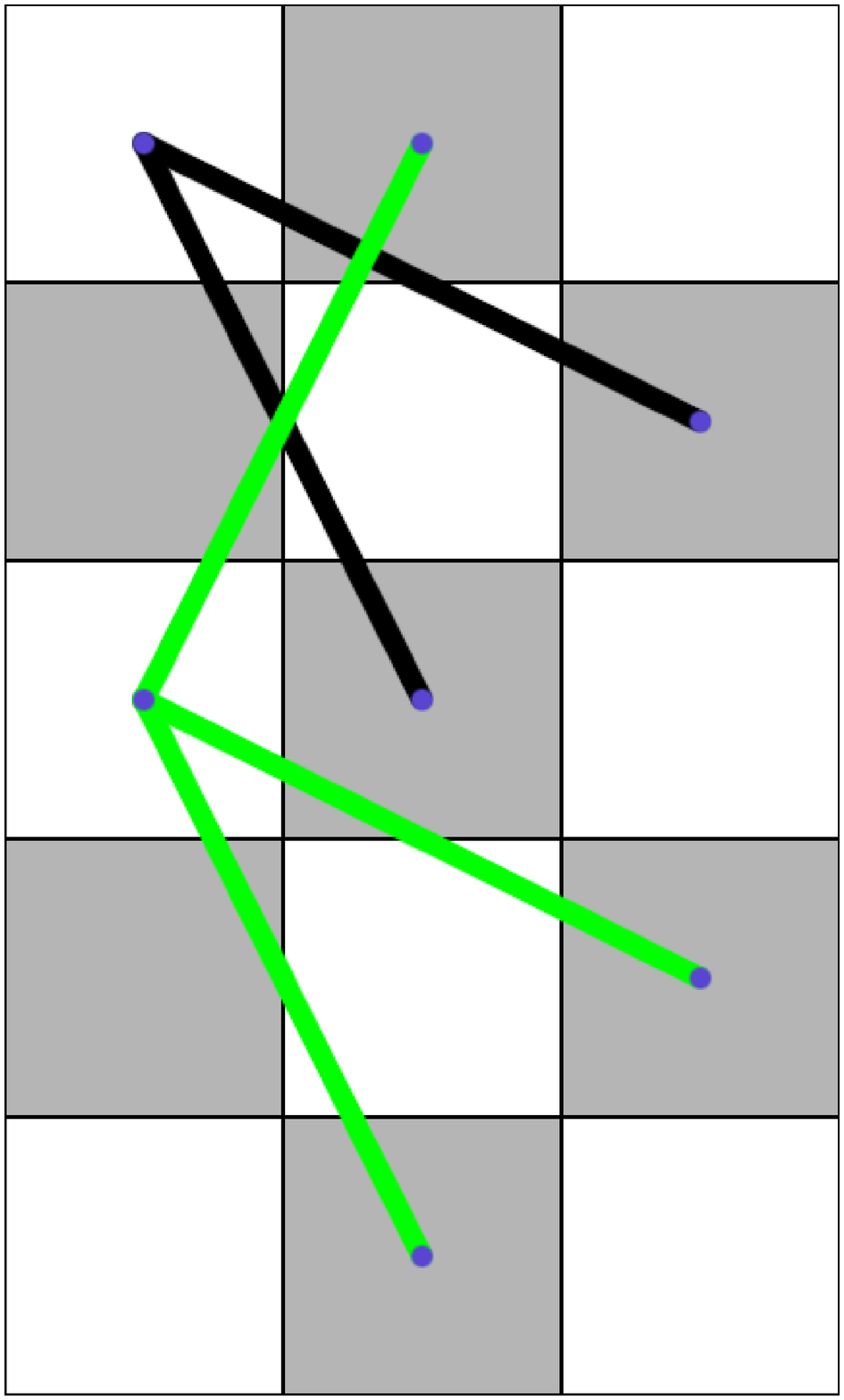}
\end{center}
\ \\
By a similar argument a site exists in each corner of the board and so every $n \times m$ tour is bi-sited, the result then follows by applying Theorem 4. \qed
\ \\
As an illustration here is an example of a bi-sited $3 \times 10$ where the sites are the highlighted edges.
\begin{center}
\includegraphics[scale=0.2]{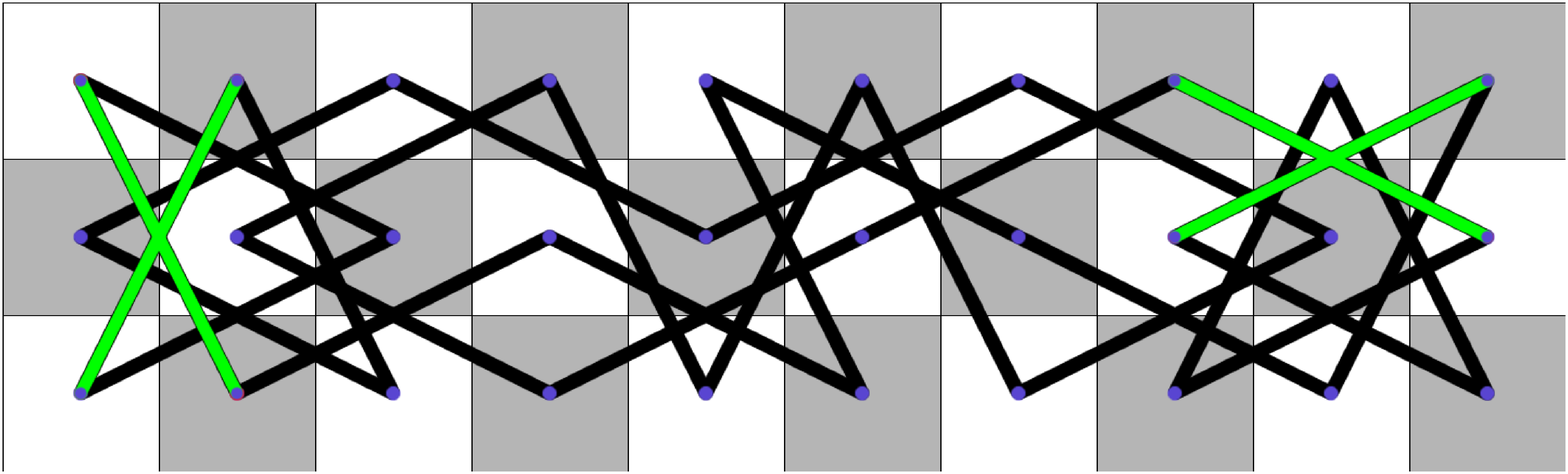}
\end{center}
\ \\
So we aim to classify all tourable chessboards by constructing bi-sited examples in small enough dimensions. It has been shown by DeMaio and Mathew ${\cite{D2007,DM2011}}$ that in $3$ dimensions the only chessboards that do not admit a knight's tour are $p \times 2 \times 2$, $p \times q \times r$ for $p,q,r \leq 3$ or $p,q,r$ all odd.\\
\ \\
{\bf Lemma 6 : }\\
A bi-sited $p \times q \times r$ $(p \geq q \geq r)$ tour exists if the following conditions hold:\\
1) $p,q$ or $r$ is even;\\
2) $p \geq 4$;\\
3) $q \geq 3$.\\
\ \\
{\bf Proof :}\\
If the first condition does not hold then no tour exists by a simple parity consideration, if any of the other conditions do not hold then $K(p,q,r)$ is disconnected. It is a simple, albeit lengthy, check that the tours constructed in ${\cite{DM2011}}$ are all bi-sited. For completeness constructions of these tours can also be found in the Appendix. \qed
\ \\
{\bf Theorem 3 :}\\
For $r \geq 3$ an $n_1 \times n_2 .... \times n_r$ $(n_1 \geq n_2 ... \geq n_r)$ tour exists if and only if the following conditions hold:\\
1) There exists an $i$ such that $n_i$ is even;\\
2) $n_1 \geq 4$;\\
3) $n_2 \geq 3$.\\
\ \\
{\bf Proof :}\\
As above if the first condition does not hold then no tour exists by a simple parity consideration, if any of the other conditions do not hold then $K(n_1,...n_r)$ is disconnected.\\
\ \\
Given an $n_1 \times n_2 .... \times n_r$ chessboard such that some $n_i$ is even then, unless (after re-ordering) $n_i = 2$ for all $i >1$ or $n_i \leq 3$ for all $i$, there is some triple, wlog $n_1, n_2, n_3$, such that $n_1 \geq 4\,,\,n_2\geq 3$ and one of $n_1 , n_2 , n_3$ is even.\\
\ \\
Therefore, by Lemma 6, a bi-sited $n_1 \times n_2 \times n_3$ tour exists and hence, by Theorem 4, a bi-sited $n_1 \times n_2 .... \times n_r$ tour does. \qed
\ \\
An immediate consequence of this is\\
\ \\
{\bf Corollary 7 :}\\
For $r\geq 3$, an $n_1 \times n_2 .... \times n_r$ tour exists if and only if $K(n_1,n_2...n_r)$ is connected.\\
\ \\
\ \\
\ \\
\section{Generalised knight's tours on a chessboard}
The knight's tour is a specific case of many general questions. A natural one to ask would be, what about move general moves? For example instead of the knight being able to move $(\pm 1, \pm 2)$ or $(\pm 2, \pm 1)$ what if the knight could move $(\pm a, \pm b)$ or $(\pm b, \pm a)$?\\
\ \\
We define an \emph{$(a,b)$-tour} of an $n_1 \times n_2...\times n_r$ chessboard to be a closed tour of the board only using moves of the form $(\pm a, \pm b)$ or $(\pm b, \pm a)$, and $K_{a,b}(n,m)$ in the obvious way. Similarly we define an \emph{$a$-site} to be a pair of lines in an $(a,b)$-tour whose endpoints are both $a$ squares away from each other, and a \emph{$b$-site} in the same way. We first look at the case where $a=1$ for ease of presentation. Note that even in this case it is not known, except when $b=2$, for which $n,m$ an $n \times m$ $(1,b)$-tour exists.\\
\ \\
{\bf Theorem 8 :}\\
If an $n \times m$ $(1,b)$-tour exists with $n > 2b + 1$ then a $n \times m \times p_1 \times p_2 .... \times p_r$ $(1,b)$-tour exists for any $r$ and any $p_1,...,p_r$.\\
\ \\
{\bf Proof :}\\
As in the proof of Corollary 5 consider the vertex $(0,0)$ in an $(1,b)$-tour. If such a tour exists, the lines $\big( (0,0) , (1,b) \big)$ and $\big( (0,0) , (b,1) \big)$ must be included, since the vertex $(0,0)$ has degree 2 in $K_{1,b}(n,m)$. Furthermore of the four lines adjacent to the point $(0,b)$ (which end at $(1,0)$, $(b,b-1)$, $(b,b+1)$ and $(1, 2b)$) at least $2$ must be in the tour, but $3$ of them form sites with the two forced lines.\\
\ \\
So, as long as the chessboard is sufficiently large to ensure that the $b$-sites in each corner are disjoint, any tour must contain at least two disjoint $b$-sites.\\
\ \\
Hence by the same argument as in Section $2$ we can construct $n \times m \times p_1 ... \times p_r$ $(1,b)$-tours for all $p_i, r \in \bb{N}$. \qed
\ \\
{\bf Theorem 9 :}\\
If $n \times m$ $(a,b)$-tours exists for all sufficiently large $n,m$ (with $n$ even) then they also exist for all $n_1 \times n_2  .... \times n_r$ for sufficiently large $n_i$ (with $n_1$ even).\\
\ \\
{\bf Proof :}\\
By a similar argument to the $(1,b)$ case, for sufficiently large $n,m$ both an $a$-site and a $b$-site must exist in all four corners of an $n \times m$ $(a,b)$-tour. In each corner these two sites are not necessarily edge disjoint but we call the union of them an \emph{$(a,b)$-site}.\\
\ \\
We claim that if an $n_1 \times n_2 ... ,\times n_s$ $(a,b)$-tour exists with $4$ $(a,b)$-sites then an $n_1 \times n_2 ... ,\times n_s \times p$ $(a,b)$-tour exists with $4$ $(a,b)$-sites for all sufficiently large $p \in \bb{N}$ and the result will follow by induction. By the above remark if an $n_1 \times n_2$ tour exists, for large enough $n_1 , n_2$ it will have $4$ $(a,b)$-sites. Now, given such an $n_1 \times n_2 ... ,\times n_{s-1}$ $(a,b)$-tour and $n_s \in \bb{N}$.\\
\ \\
In the case where $a=1$ we would start, as in the proof of Theorem 4, with $n_s$ cycles stacked on top of each other, where the cycles are copies of the $n_1 \times n_2 ... ,\times n_{s-1}$ $(a,b)$-tour, and at any point we could use $b$-sites to join two cycles that differed by one layer in the stack into a longer cycle, repeating this until there is one cycle left. In this case the strategy is obvious, simply go through the layers in order.\\
\ \\
The situation for general $a$ and $b$ is similar, except we can only join cycles that are $a$ or $b$ layers apart. So instead of just going through the layers one by one we need to find a path through the stack only ever moving $a$ or $b$ layers at a time, although in fact we can get by with slightly less. Looking back at our argument in section $2$ we only ever used two sites on each layer, since we have 4 $(a,b)$-sites we only need to find a tree through the layers with maximum degree $4$.\\
\ \\
So, since for a tour to exist $a$ and $b$ must be coprime, we can use the $a$-sites to adapt our $n_s$ copies of the $n_1 \times n_2  .... \times n_{s-1}$ $(a,b)$-tour to form $b$ cycles, each touring the residue classes of layers modulo $b$ and then use the $b$-sites to join each cycle in turn. This is possible if $n_s \geq a+b+1$. \\
\ \\
This forms a tree through the layers and, since a tree has at least two leaves, that is two layers on which only $1$ $(a,b)$-site is used, there are still at least $4$ $(a,b)$-sites in the new tour.\qed
\ \\
It is not known in general for which $a,b$ $(a,b)$-tours exist on sufficiently large chessboards. Knuth ${\cite{K1994}}$ showed that in two dimensions if gcd$(a+b,a-b)=1$ then $K_{a,b}(n,m)$ is connected for $n \geq 2b$ and $m \geq a+b$ (and that this is tight), otherwise $K_{a,b}(n,m)$ is disconnected. In light of the conditional nature of Theorem 9 it seems natural to conjecture\\
\ \\
{\bf Conjecture 1 :}\\
For all $a,b$ such that gcd$(a+b,a-b)=1$ there exists $M$ such that if $n$ is even and $n,m \geq M$ then an $n \times m$ $(a,b)$-tour exists.\qed
\ \\
Even more than this it might be true that, as in the $(1,2)$ case, for sufficiently large dimensions it is enough that the $(a,b)$-knight's graph is connected to ensure a tour exists.\\
\ \\
{\bf Conjecture 2 :}\\
For all $a,b$ such that gcd$(a+b,a-b)=1$ there exists $r$ such that if $K_{a,b}(n_1,n_2...,n_r)$ is connected then an $n_1 \times n_2 ... \times n_r$ $(a,b)$-tour exists.\qed
\ \\
Finally we can consider even more general knight's move. For example given an $s$-tuple $a_1,...a_s$ we can consider $(a_1,...,a_s)$-tours on $n_1 \times n_2... \times n_r$ chessboards as long as $r \geq s$. It is not hard to show that given $a_1...,a_s$ and $r > s$ $K_{a_1...,a_s}(n_1...,n_r)$ is connected, for sufficiently large $n_i$, if and only if $\sum a_i \equiv 1$ mod$(2)$ and gcd$(a_1...,a_s)=1$. If $r=s$ then we require the additional constraint that at least one of the $a_i$ are even. With this in mind we conjecture\\
\ \\
{\bf Conjecture 3 :}\\
For all $a_1,...a_s \in \bb{N}$ such that $gcd(a_1,a_2...a_s)=1$ and $\sum a_i \equiv 1$ mod$(2)$ there exists $M$ such that if $n_1,n_2,...n_r \geq M$, with $r > s$, and $n_1$ is even then an $n_1 \times n_2 ... \times n_r$ $(a_1,a_2...a_s)$-tour exists.\qed
\newpage
\section*{Appendix}
{\bf Lemma 6 : }\\
A bi-sited $p \times q \times r$ $(p \geq q \geq r)$ tour exists if the following conditions hold:\\
1) One of $p,q$ or $r$ is even;\\
2) $p \geq 4$;\\
3) $q \geq 3$.\\
\ \\
{\bf Proof :}\\
We will start by constructing all possible $2$-dimensional tours since, by using Corollary 5, this will provide constructions of a large class of bi-sited $3$-dimensional tours.\\
\ \\
We call a tour \emph{seeded} if it includes the edges $\big( (0,m-3)\,,\,(1,m-1) \big)$ and $\big( (n-3,0)\,,\,(n-1,1) \big)$. For example the the $10 \times 3$ tour below is seeded:\\
\begin{center}
\includegraphics[scale=0.2]{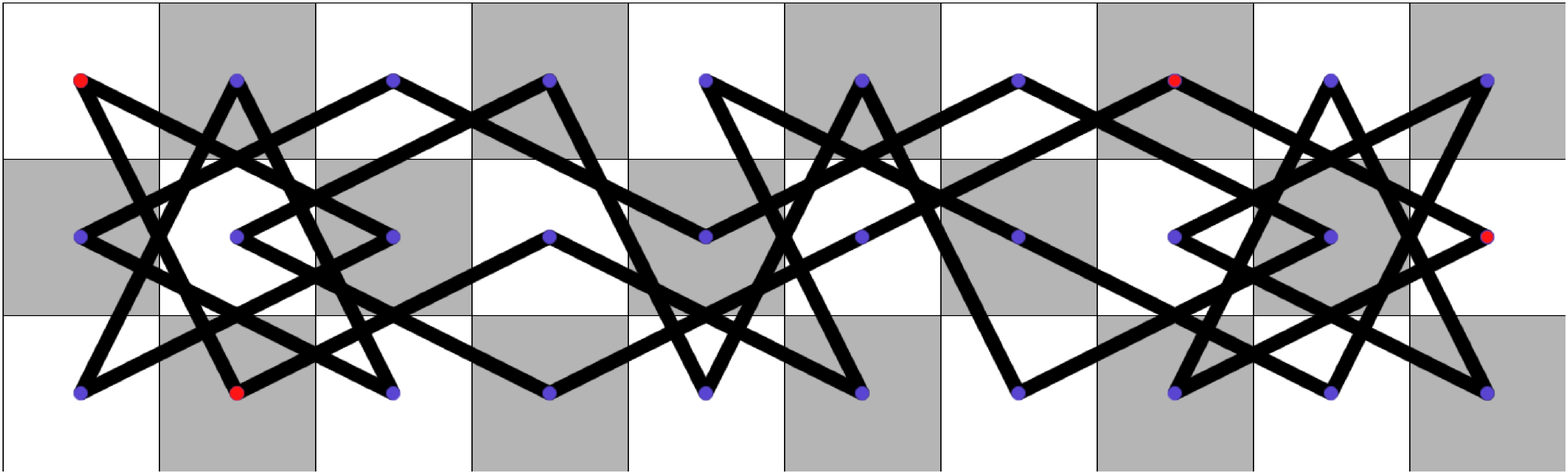}
\end{center}
\ \\
We call an open tour on a $4 \times m$ chessboard a $4 \times m$ \emph{extender} if it is a tour starting at $(3,m-1)$ and ending at $(3,m-2)$, like the $4 \times 3$ extender below:\\
\begin{center}
\includegraphics[scale=0.2]{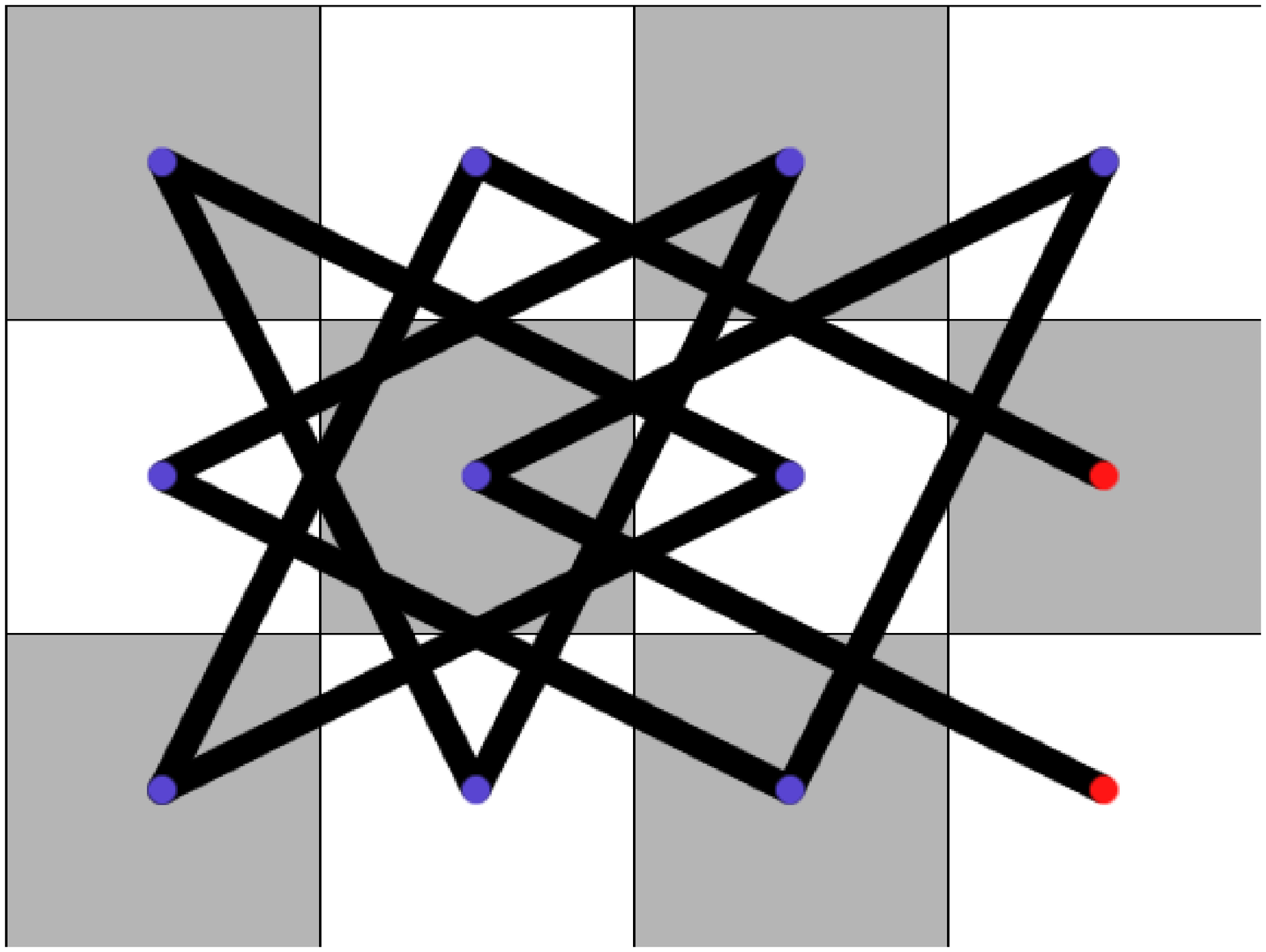}
\end{center}
\ \\
{\bf Lemma 10 :}\\
There exists a seeded $4 \times m$ extender for all  $m \neq 1,2$ or $4$.\\
\ \\
{\bf Proof :}\\
Observe that if we place the following pattern:
\begin{center}
\includegraphics[scale=0.2]{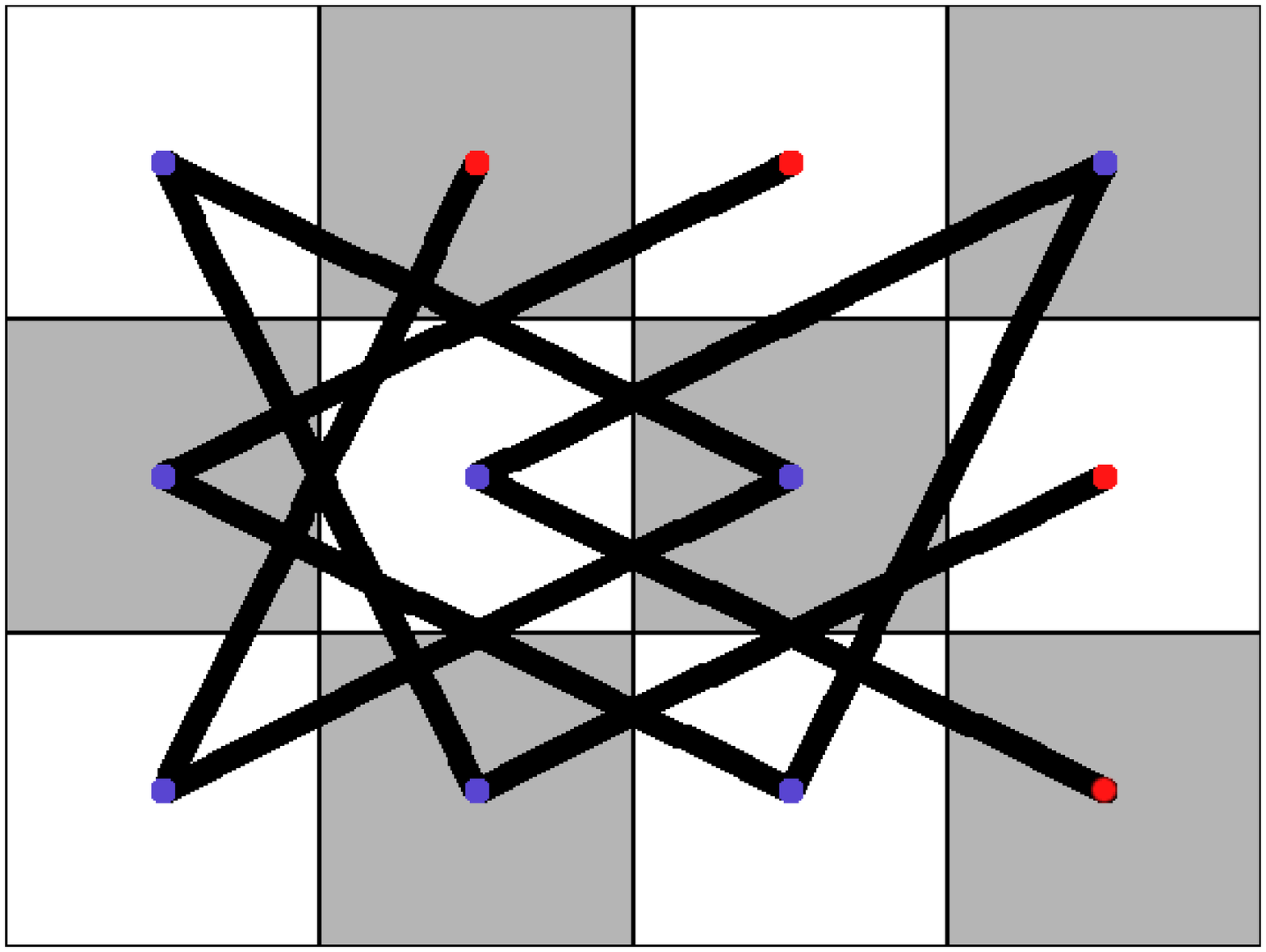}
\end{center}
\ \\
below a seeded $4 \times m$ extender and add the lines $\big( (3,m-2)\,,\, (2,m) \big)$ and $\big( (3,m-1)\,,\, (1,m) \big)$ then it will form a seeded $4 \times (m+3)$ extender. So, along with the seeded $4 \times 3$ extender above, it remains to exhibit seeded $4 \times 5$ and $4 \times 7$ extenders, as below:
\begin{center}
\includegraphics[scale=0.2]{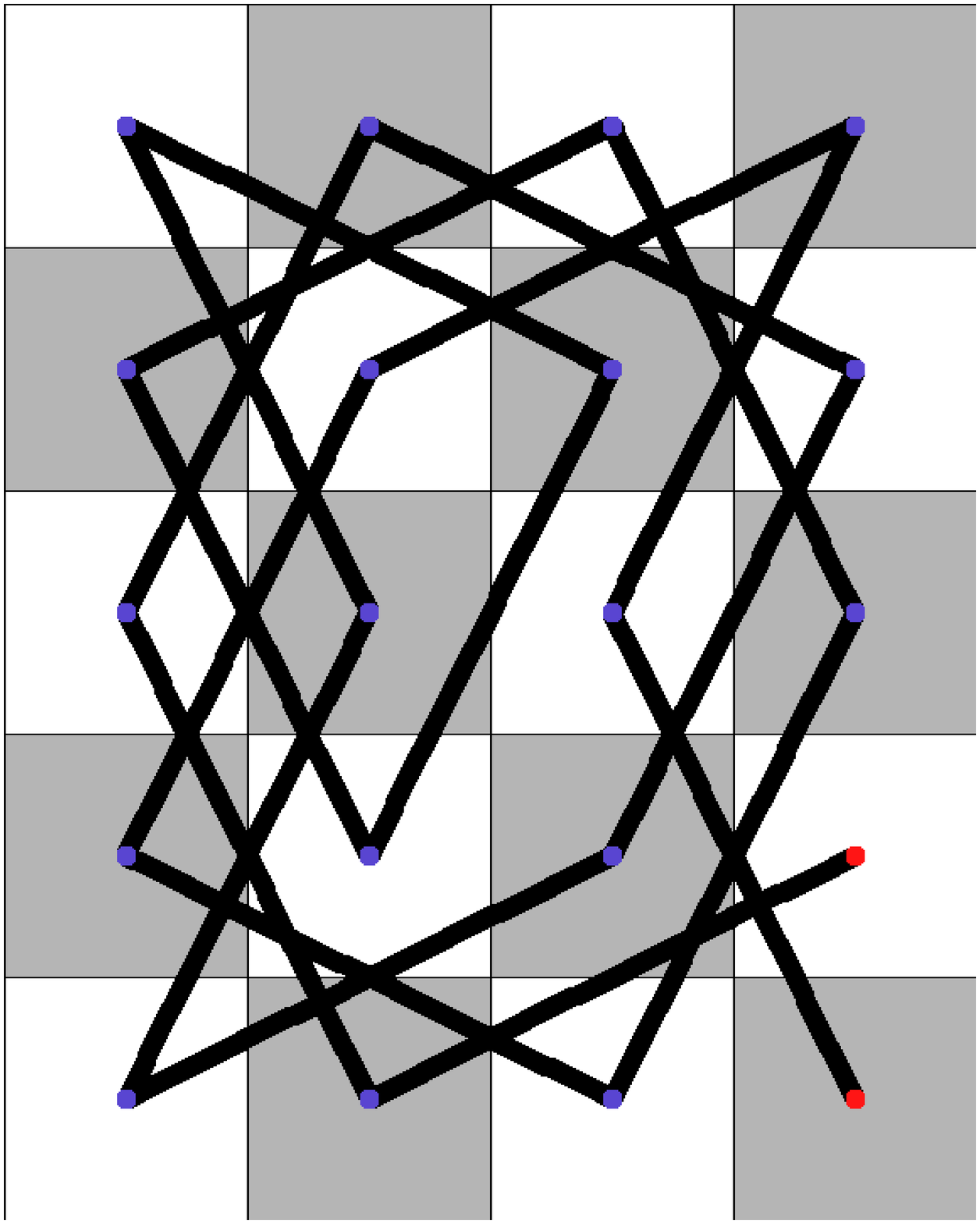} \includegraphics[scale=0.2]{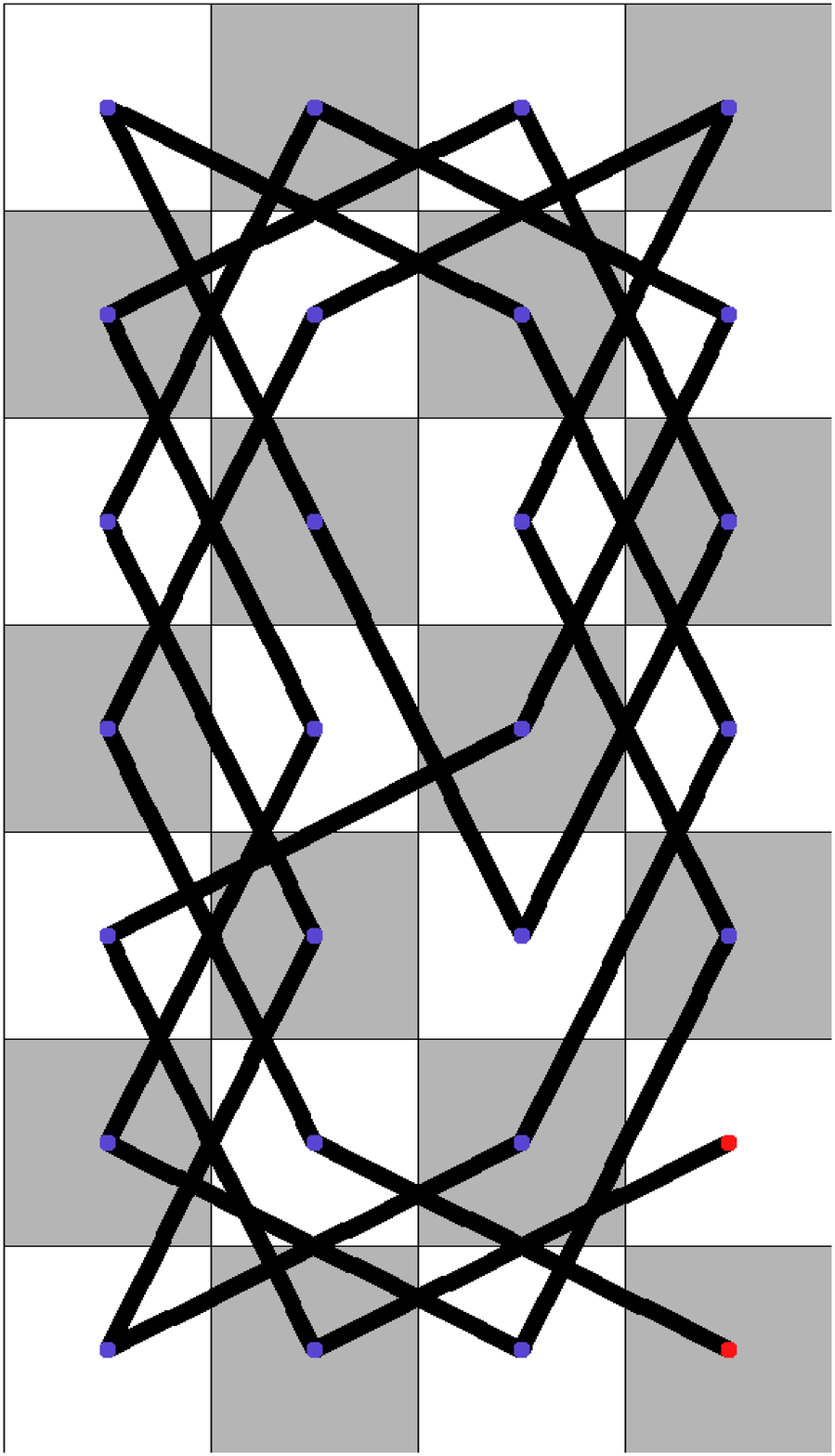}
\end{center}
\qed
\ \\
{\bf Lemma 11 :}\\
If a seeded $n \times m$ tour exists then a seeded $(n + 4) \times m$ tour exists.\\
\ \\
Note also that, since a seeded $n \times m$ tour is equivalently a seeded $m \times n$ tour after a suitable reflection, this lemma will imply that if a seeded $n \times m$ tour exists then so does a seeded $(n + 4k) \times (m + 4l)$ tour for all  $k,l \in \bb{N}$.\\
\ \\
{\bf Proof :}\\
Given a seeded $n \times m$ tour with $m \neq 1,2$ or $4$ then there exists a seeded $4 \times m$ extender. Now if we place a seeded $4 \times m$ extender to the left of a seeded $n \times m$ tour like below:
\begin{center}
\includegraphics[scale=0.2]{3x4.eps} \includegraphics[scale=0.2]{3x10.eps}
\end{center}
\ \\
By removing the line $\big( (5, m-1)\,,\,(4,m-3)\big)$ and adding in the two lines $\big( (3,m-1) \,,\,(4,m-3)\big)$ and $\big( (5, m-1)\,,\,(3,m-2)\big)$ we form a $(n + 4) \times m$ tour. Note also that this tour is still seeded. Hence a seeded $(n + 4) \times m$ tour exists. \qed
\ \\
{\bf Lemma 12 : }\\
An $n \times m$ $(n \geq m)$ tour exists if the following conditions hold:\\
1) $n$ or $m$ is even;\\
2) $m \not\in \{1,2,4\}$;\\
3) $(n,m)$ $\neq (4,3)\,,\,(6,3)$ or $(8,3)$.\\
\ \\
{\bf Proof :}
By Lemma 11 it is sufficient to exhibit a seeded $n \times m$ tour for all different pairs of residue modulo $4$ (excepting the cases where both are odd), and possibly some small cases. A quick check will show it is sufficient to use as base cases seeded $3 \times 10$, $3 \times 12$, $5 \times 6$, $5 \times 8$, $6 \times 6$, $6 \times 7$, $6 \times 8$, $7 \times 8$ and $8 \times 8$ tours, which appear below: 
\ \\
\begin{center}
\includegraphics[scale=0.2]{3x10.eps}
\includegraphics[scale=0.2]{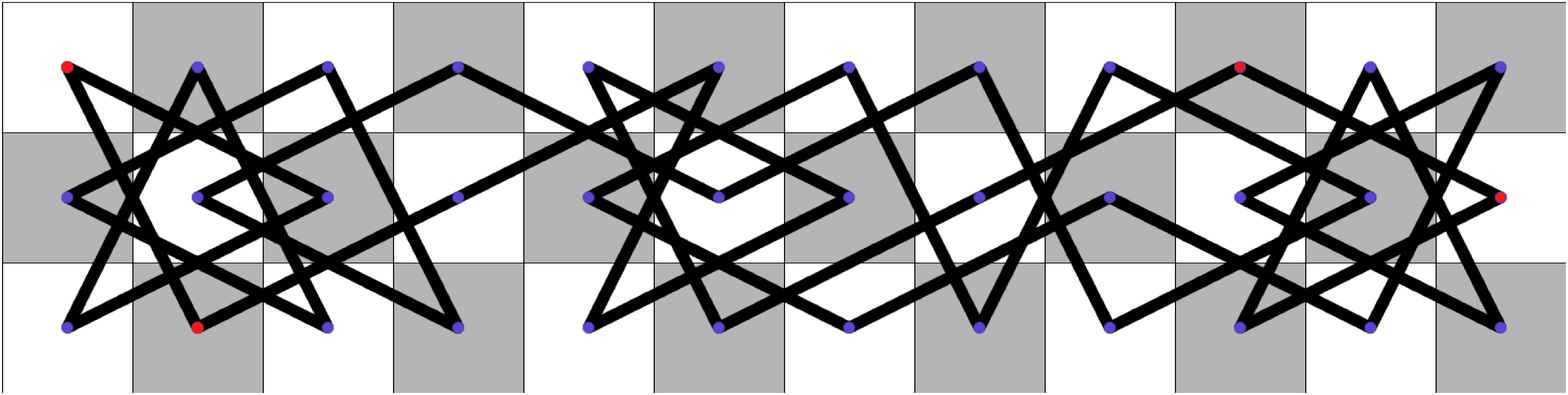}
\includegraphics[scale=0.2]{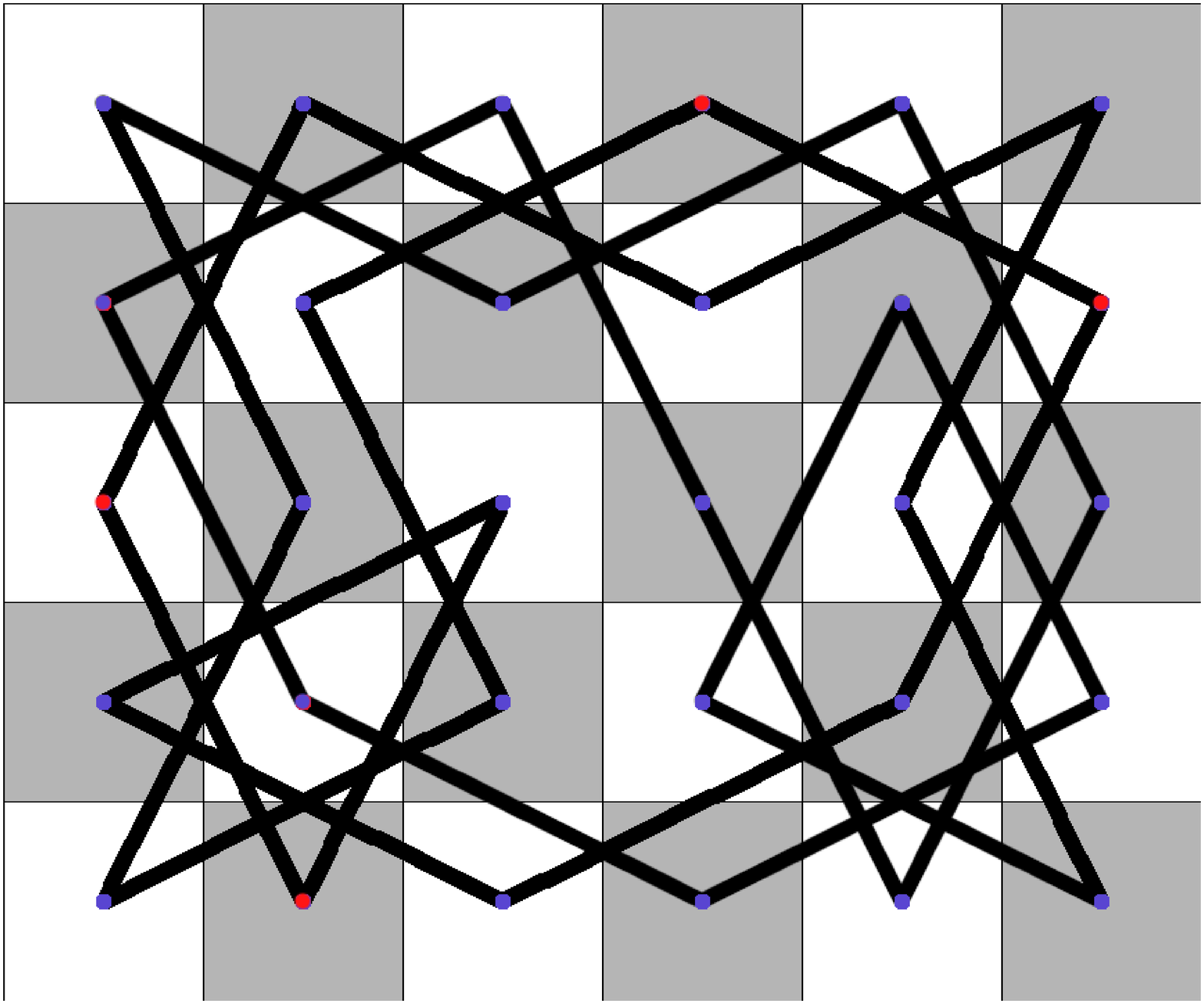}
\includegraphics[scale=0.2]{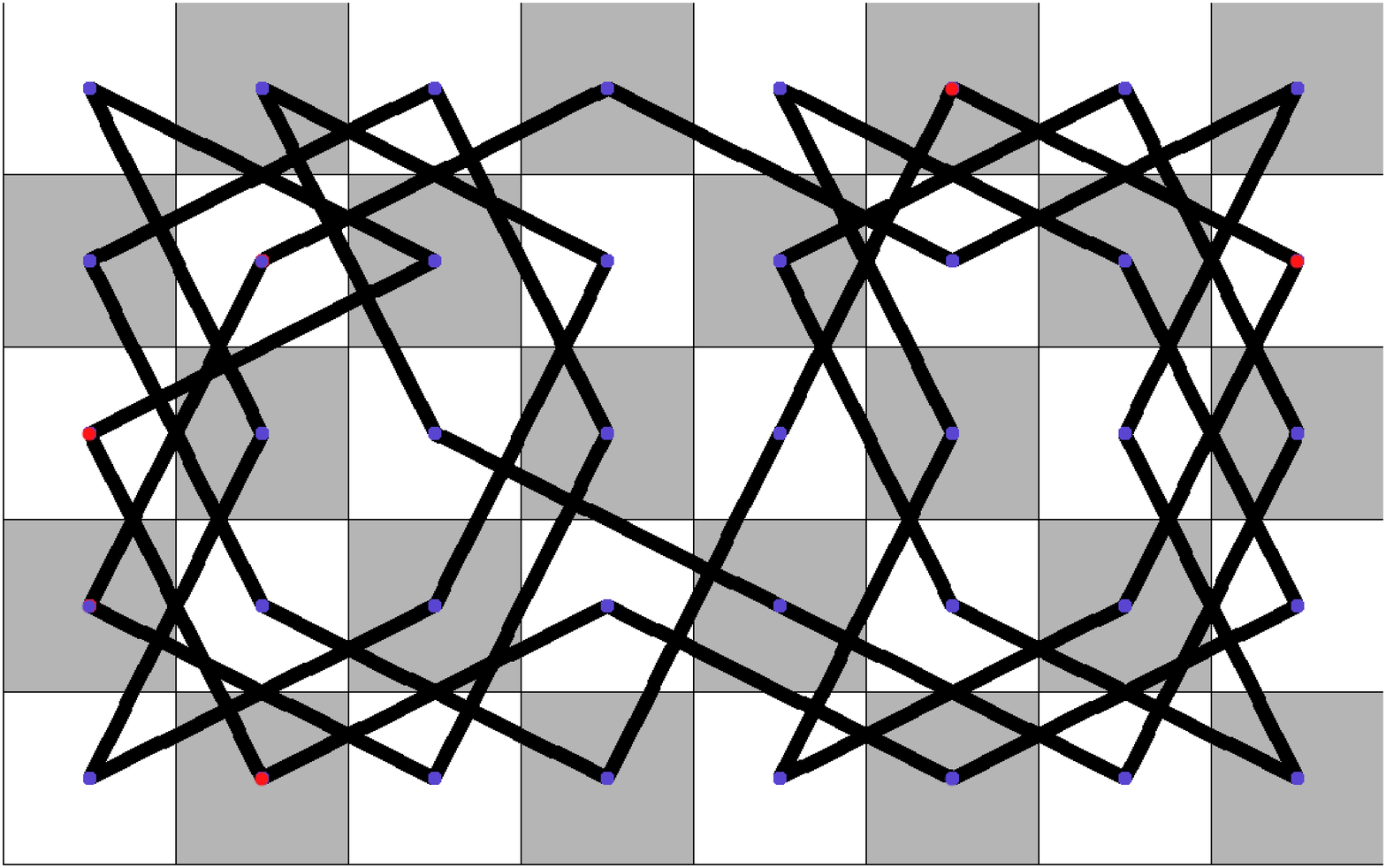}
\includegraphics[scale=0.2]{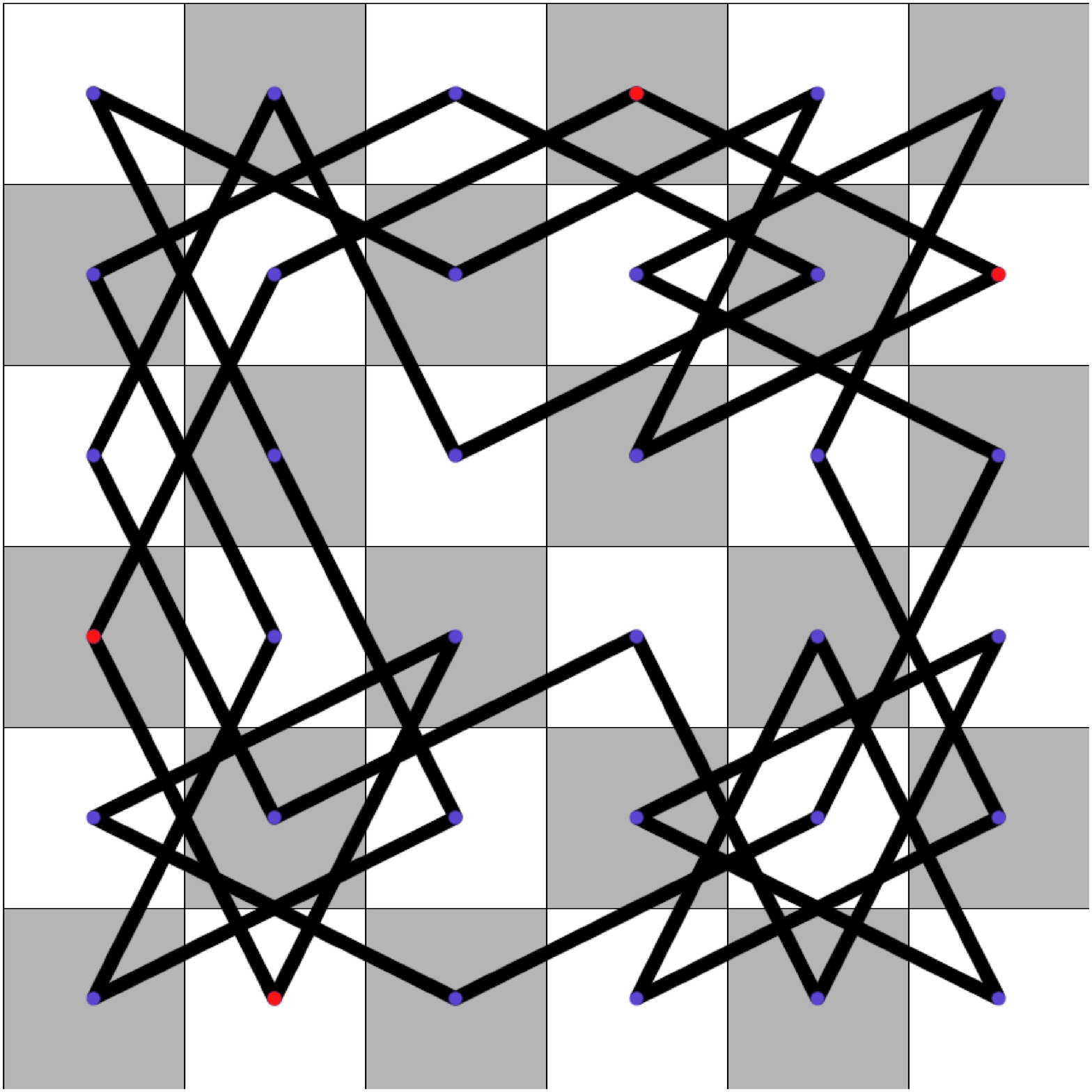}
\includegraphics[scale=0.2]{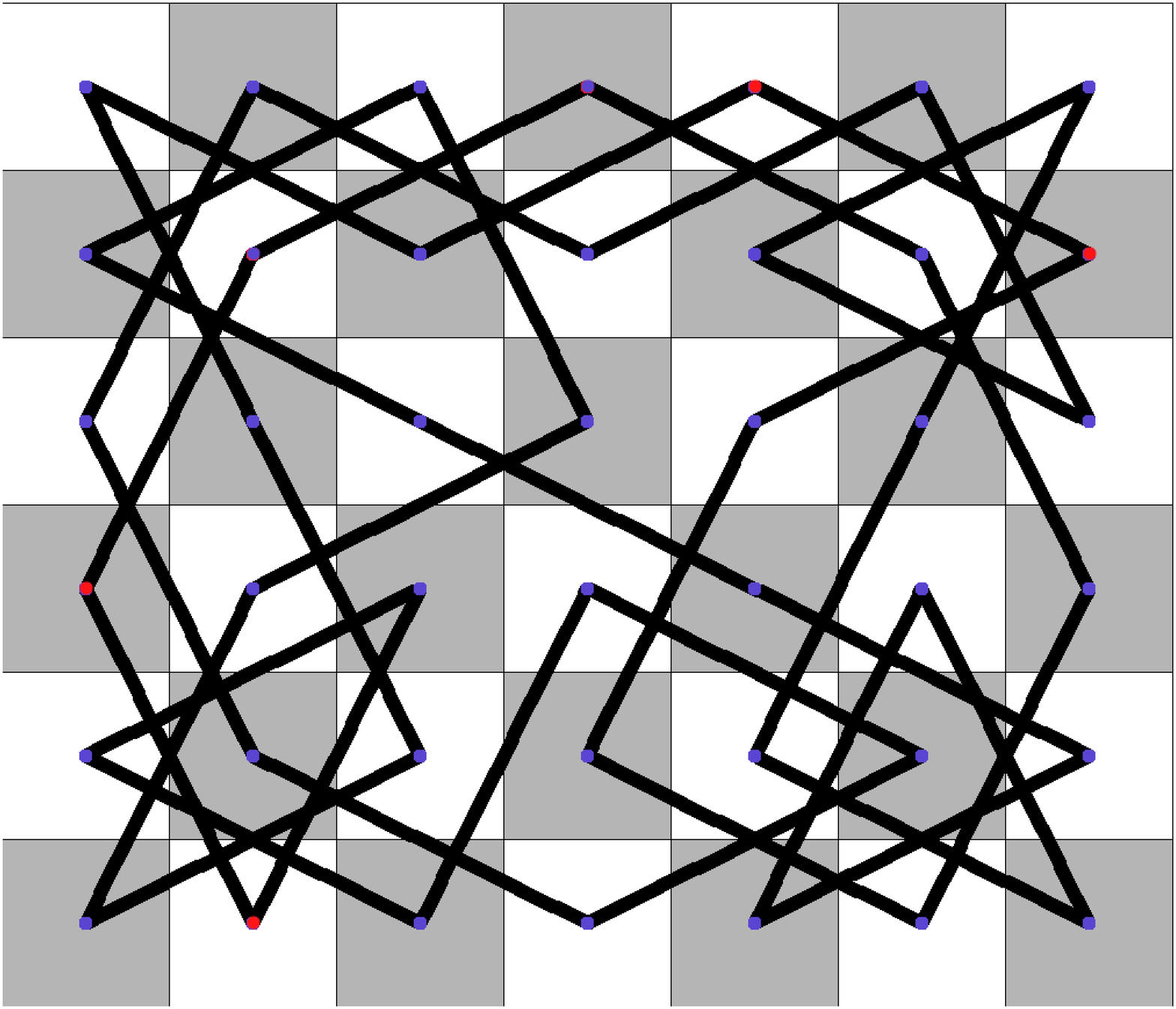}
\includegraphics[scale=0.2]{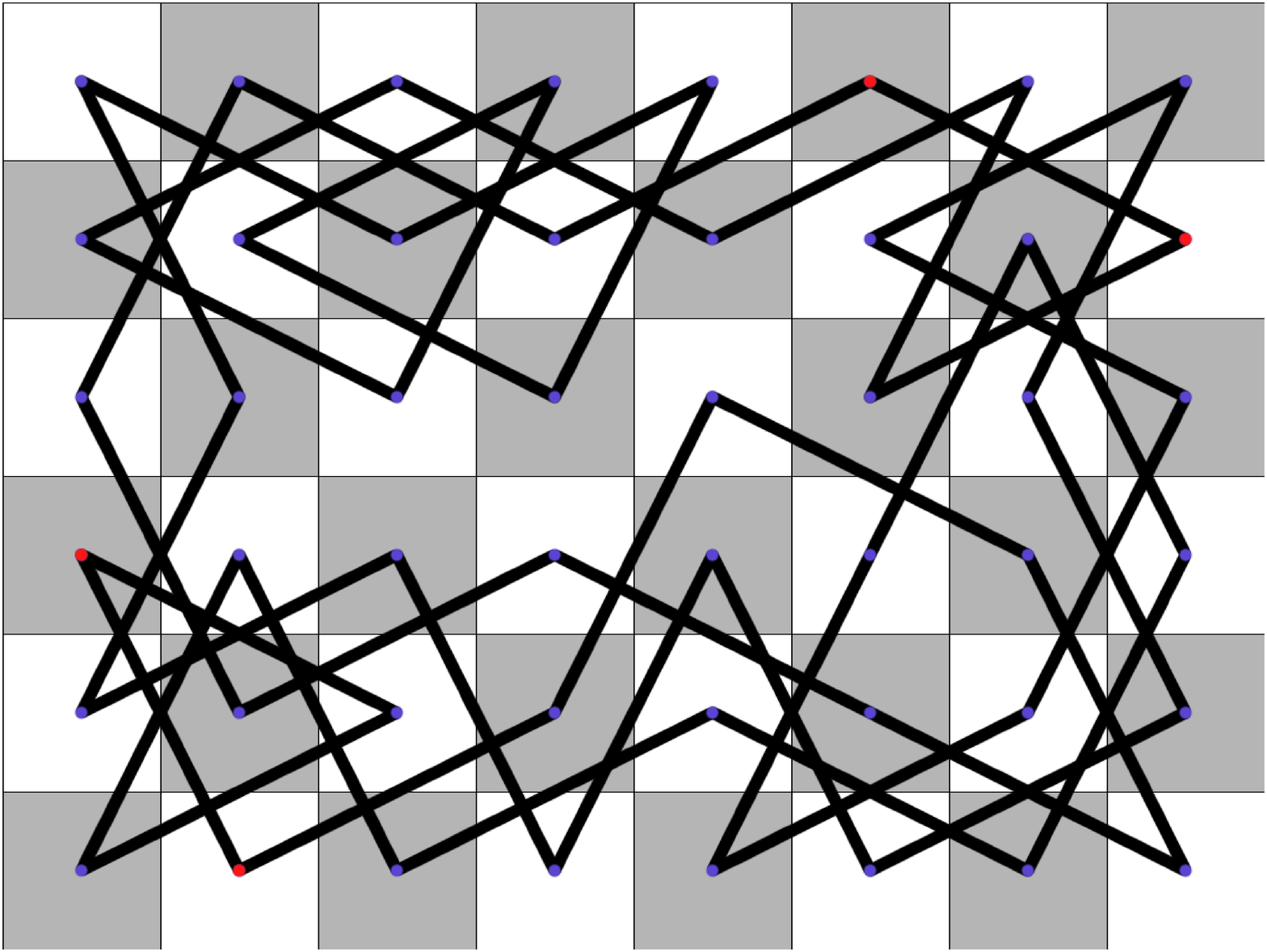}
\includegraphics[scale=0.2]{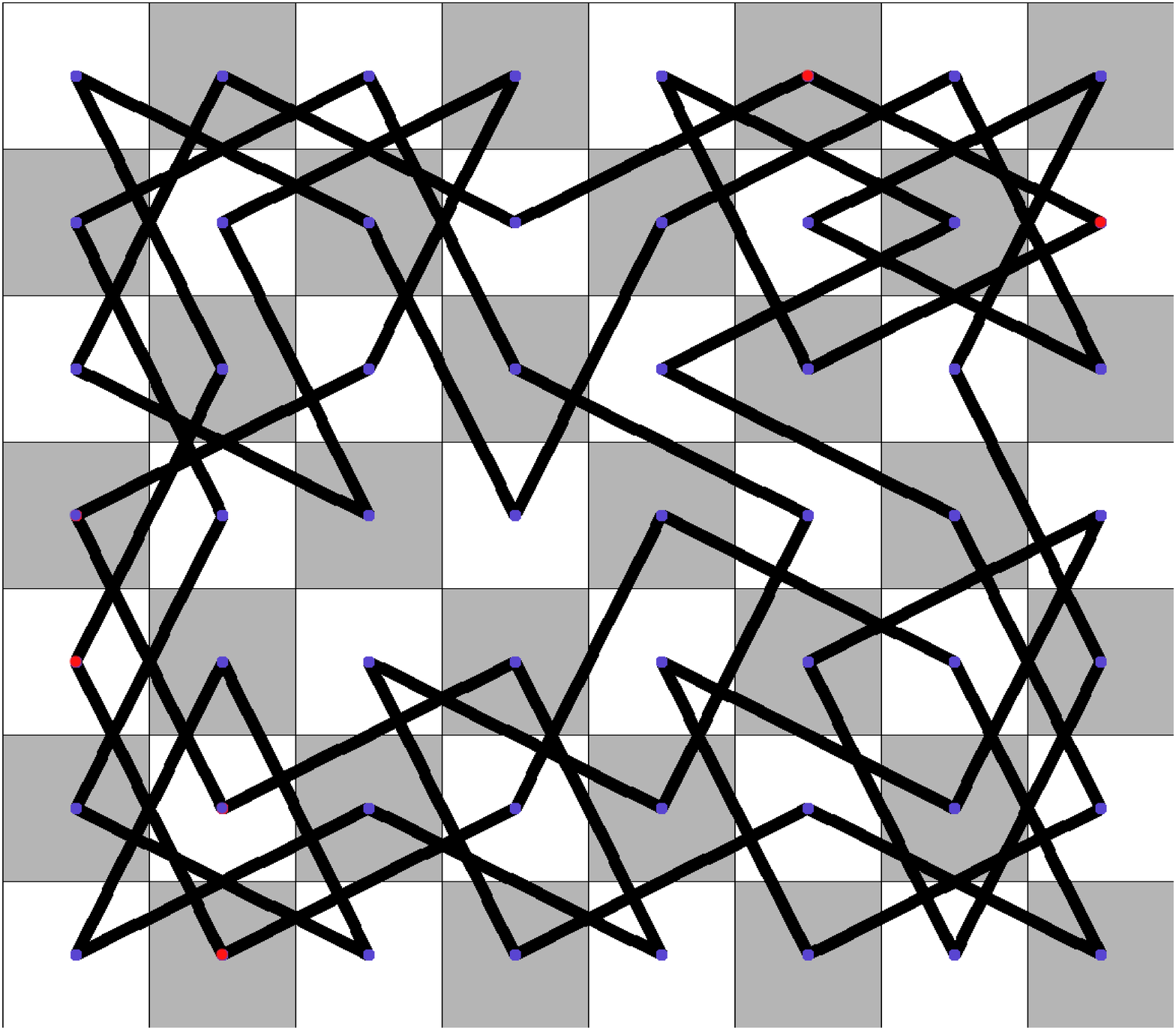}
\includegraphics[scale=0.2]{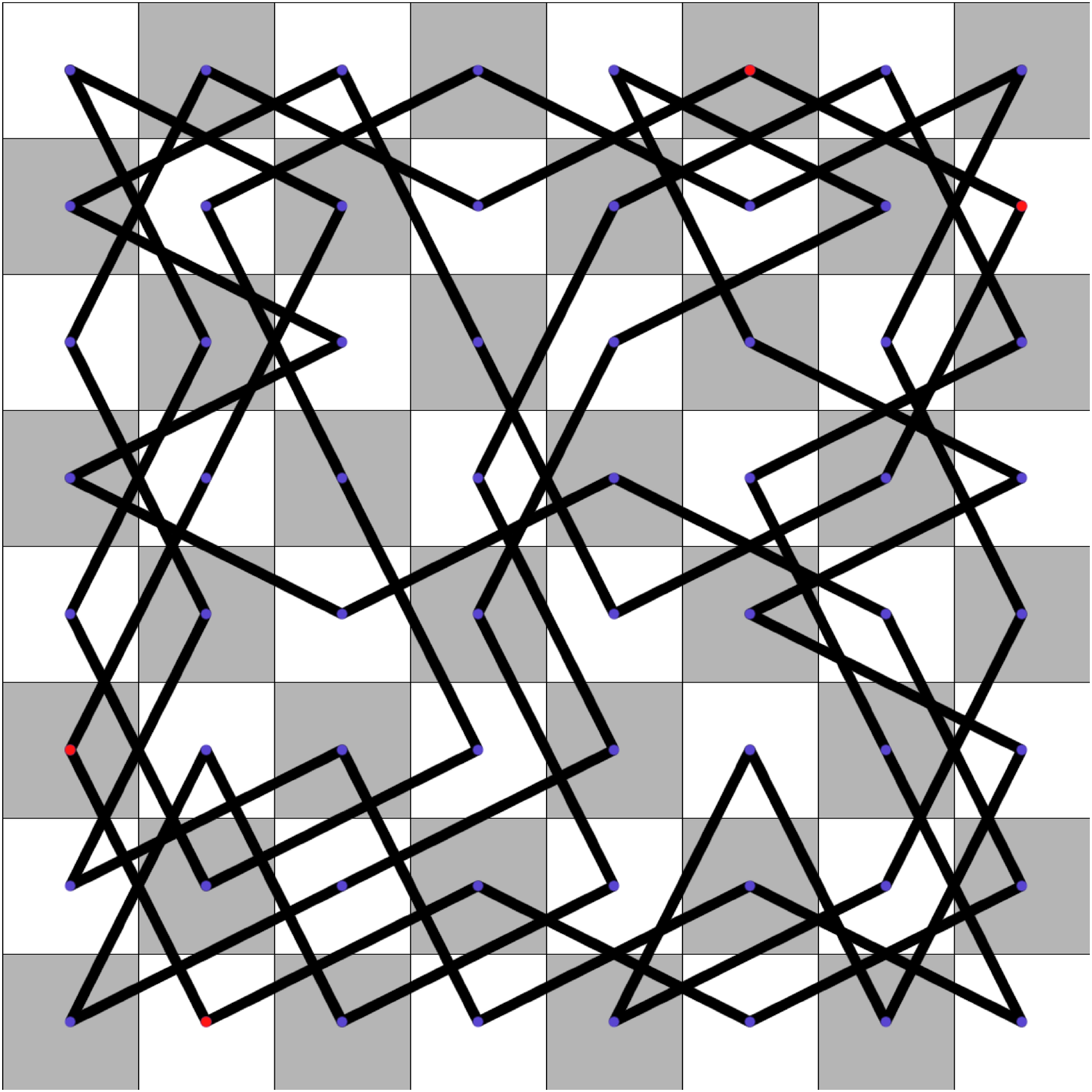}
\end{center}
\qed
\ \\
So by Corollary 5 we can construct $3$-dimensional bi-sited tours for all $n \times m \times p$ when an $n \times m$ chessboard admits a knight's tour. Also trivially if such a tour exists for $n_1 \times n_2 .... \times n_r$ then it also exists for $n_{\phi(1)} \times n_{\phi(2)} ... \times n_{\phi(r)}$ for any permutation $\phi$ of $\{1,2...,r\}$ so the order is irrelevant. We will split the remaining tours into cases, let us first consider $n \times m \times 2$ and $n \times m \times 4$ for $n,m \geq 5$ and odd.\\
\ \\
Note that given an open tour of an $n \times m$ board that starts at $(n-1,m-1)$ and ends two squares above at $(n-1,m-3)$ we can construct a closed tour of an $n \times m \times 2$ board by putting two copies of the open tour on top of each other and adding the lines $\big((n-1,m-1,0),(n-1,m-3,1)\big)$ and $\big((n-1,m-1,1),(n-1,m-3,0)\big)$, like in the example below:\\
\begin{center}
\includegraphics[scale=0.2]{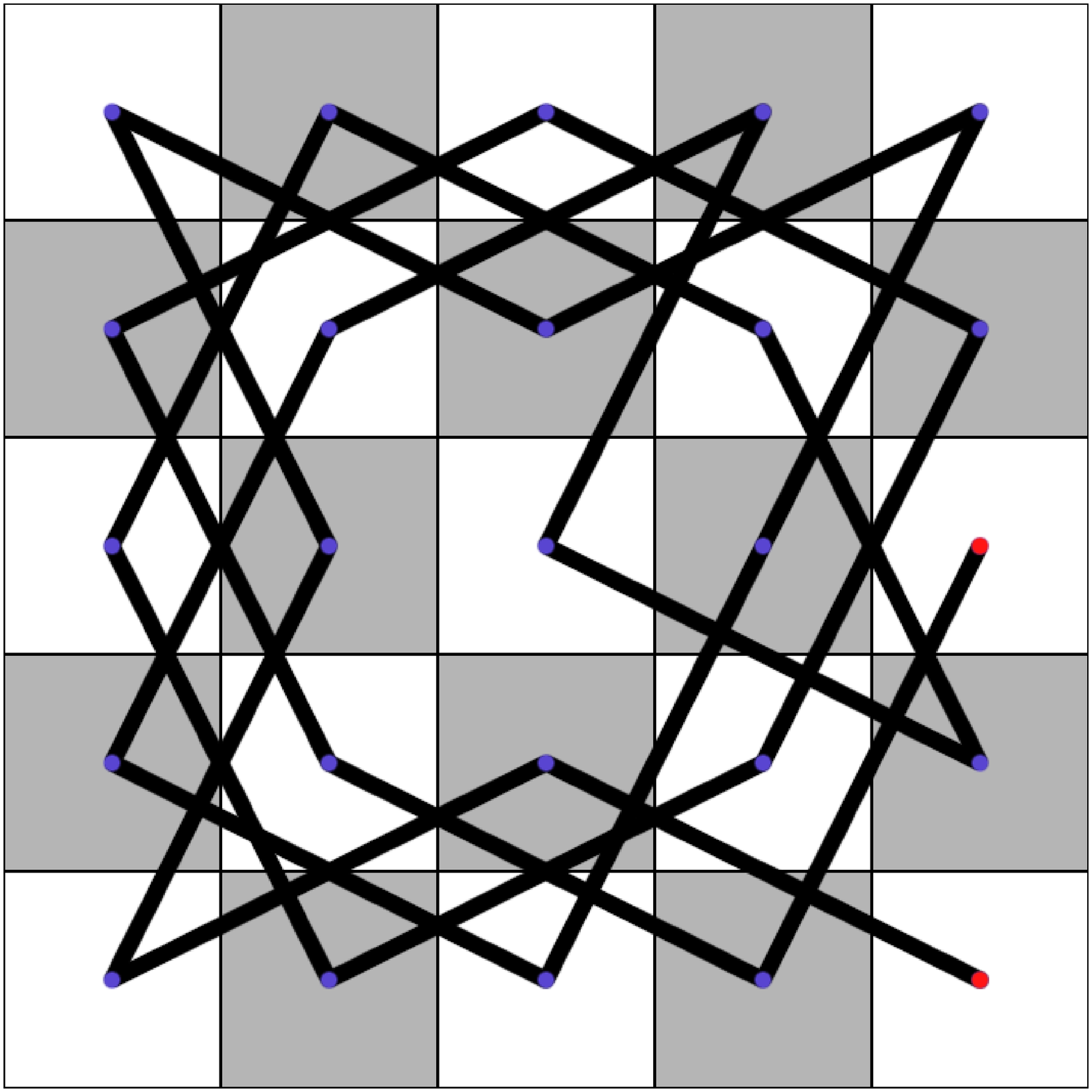}
\end{center}
\ \\
Also, using the $4 \times m$ extenders, since this tour is seeded, we could use this board to build such tours for $n \times m \times 2$ boards as long as $m,n \equiv 1$ mod$(4)$. So to complete this case we need to show that such seeded open tours exist for $5 \times 7$, $7 \times 5$ and $7 \times 7$.\\ 
\begin{center}
\includegraphics[scale=0.2]{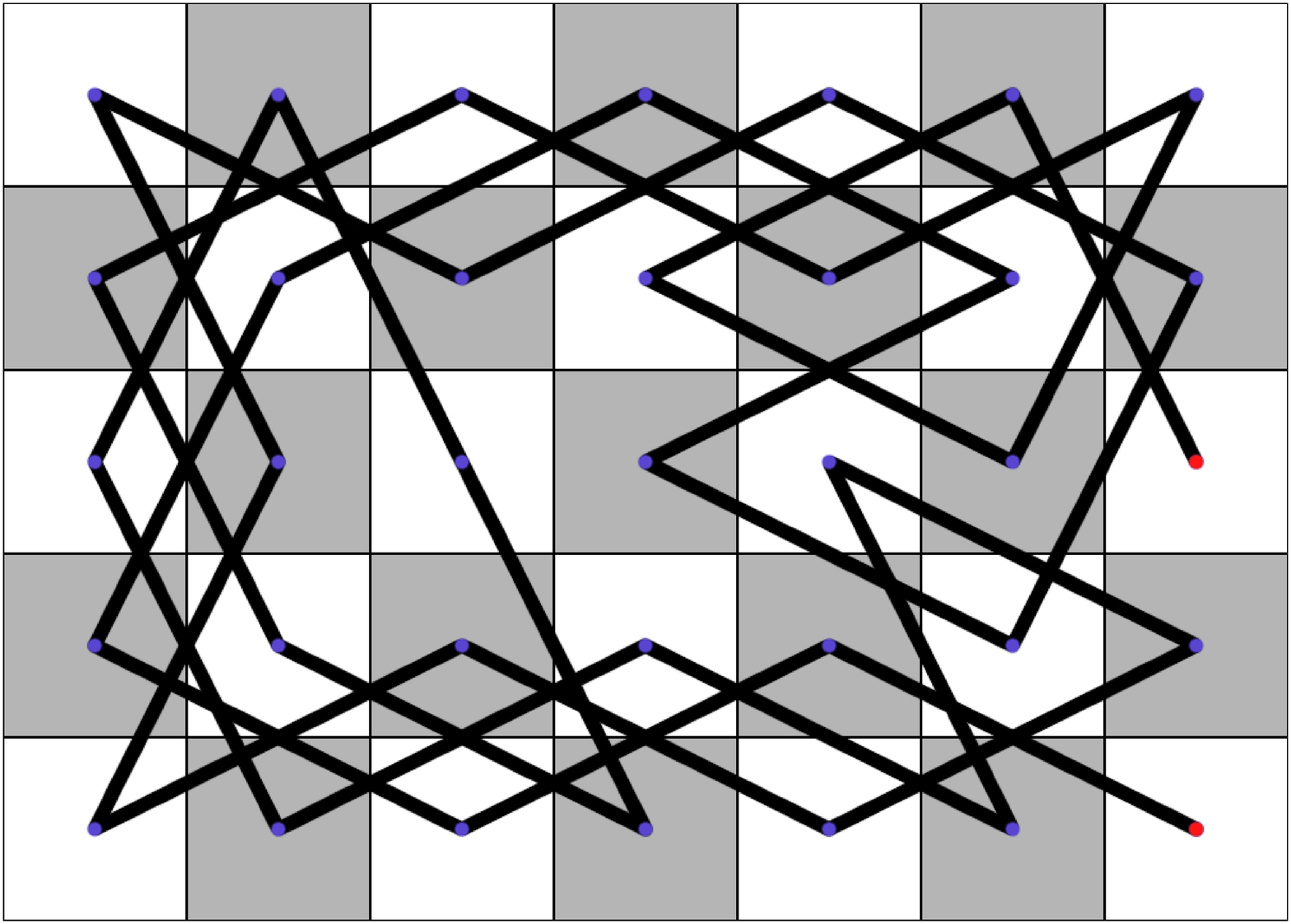}
\includegraphics[scale=0.2]{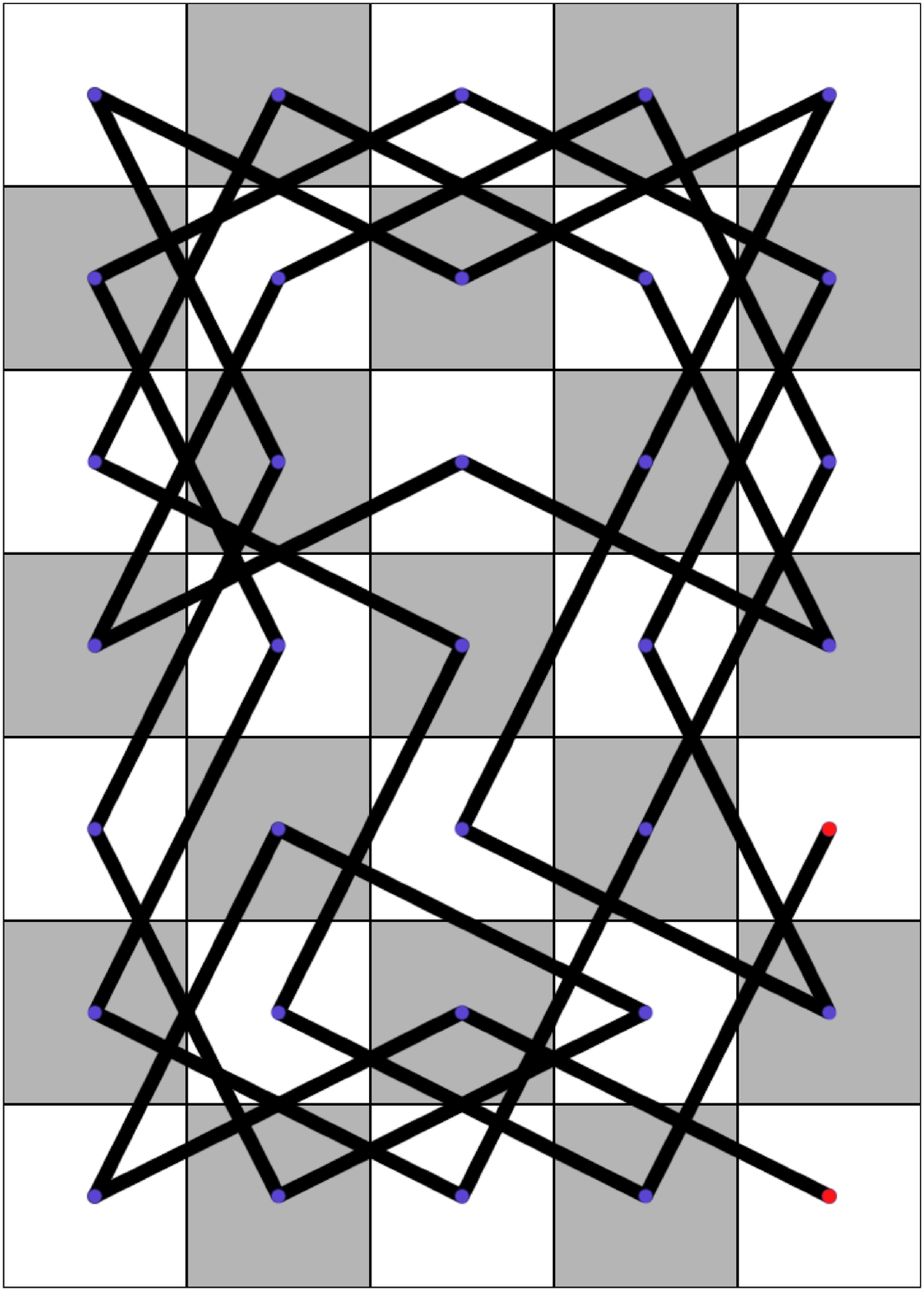}
\ \\
\ \\
\includegraphics[scale=0.2]{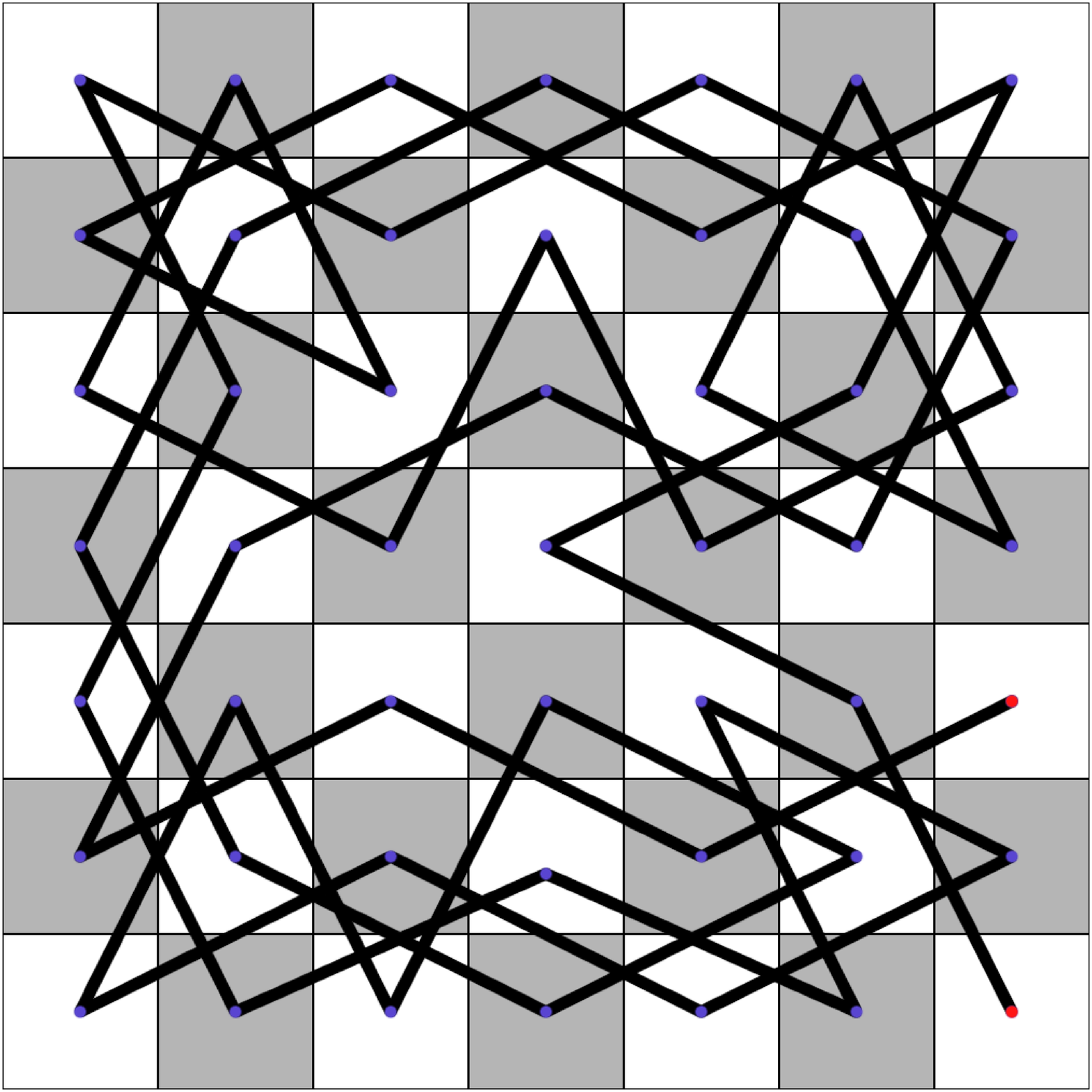}
\end{center}
\ \\
Also we need to check that these tours are bi-sited, which by construction they will be, having for example a site in the bottom left corner of both layers. Finally note that, since the two sites are directly on top of each other, using the method described in Section $2$ we can stack these to create bi-sited $n \times m \times 2k$ tours for all $n,m \geq 5$ odd and for all $k$, in particular for $n \times m \times 4$.\\
\ \\
So we are done for triples $n,m,p$ as long as $n,m \geq 5$. So all the remaining cases have at least two sides smaller than $5$. Let us consider the cases with a side of size $4$ first.\\
\ \\
The method we have used so far to draw tours will be insufficient to demonstrate more complicated $3$-dimensional tours so we will simply present them layer by layer with each square numbered with the order it appears in the tour, starting with the topmost layer to the left and so on. The colours of the squares are merely to differentiate them, not the standard black-white colouring of the chessboard.\\
\ \\
Firstly we will do $4 \times 4 \times n$. Below we exhibit a $4 \times 4 \times 2$ and a $4 \times 4 \times 3$ tour.
\begin{center}
\includegraphics[scale=0.2]{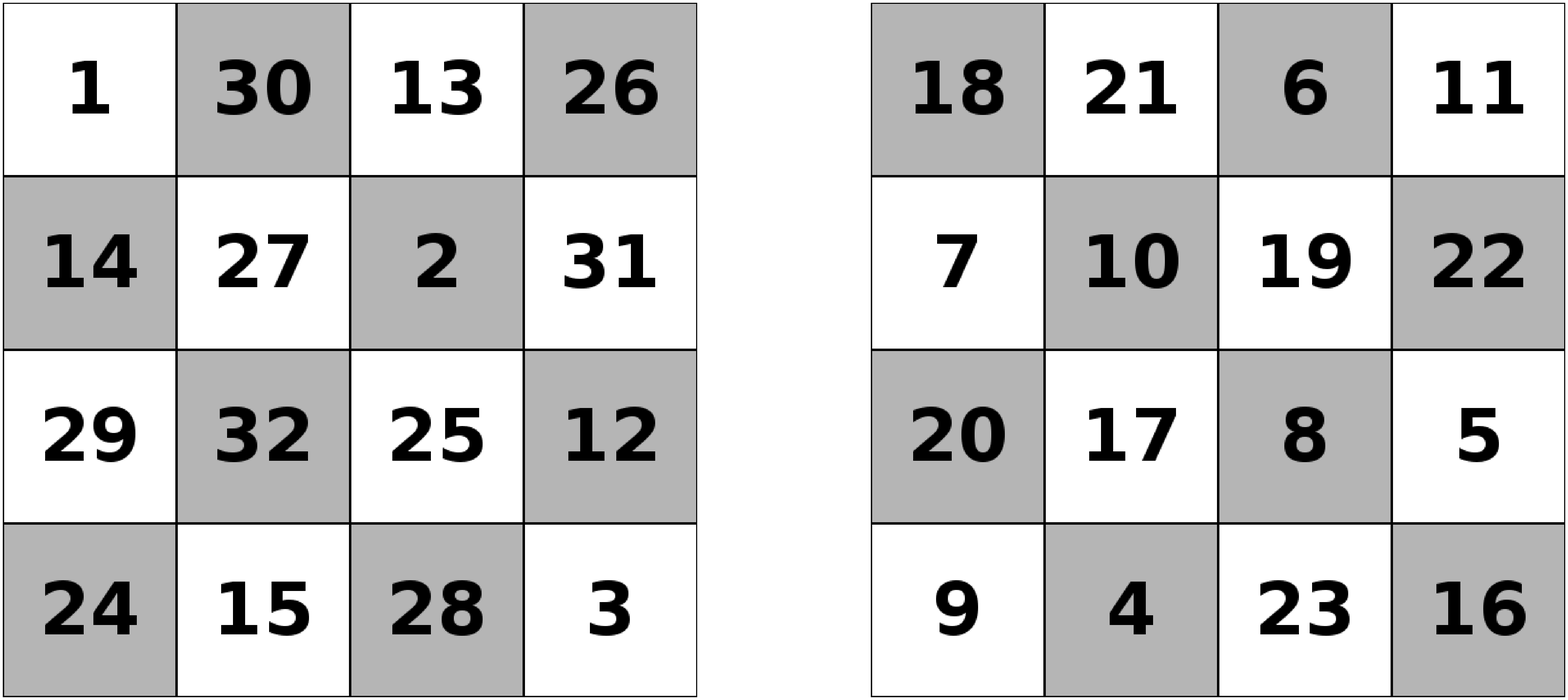}
\ \\
\ \\
\includegraphics[scale=0.2]{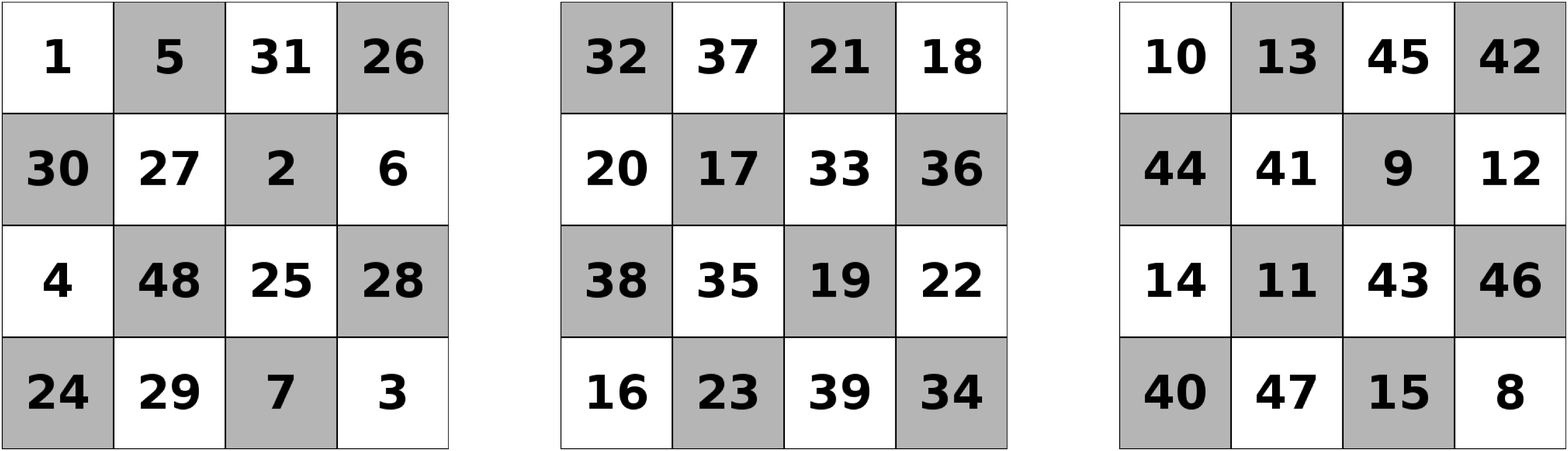}
\end{center}
\ \\
Notice the sites in the top left corners of the top and bottom layers of each of them, that is the lines $1-32$, $29-30$, $18-17$ and $20-21$ in the $4 \times 4 \times 2$ tour and the lines $1-48$, $4-5$, $10-11$ and $13-14$ in the $4 \times 4 \times 3$ tour . So we can stack any number of these on top of each other to form bi-sited $4 \times 4 \times n$ tours for all $n$.\\
\ \\
More concretely we can form a $4 \times 4 \times 4$ tour by removing the line $20-21$ from a copy of a $4 \times 4 \times 2$ tour and placing it on top of another copy with the line $1-32$ removed, then add in the lines $1-20$ and $32-21$. In a similar fashion we can add any number of $4 \times 4 \times 2$ and $4 \times 4 \times 3$ tours together.\\
\ \\
Next let us look at $4 \times 3 \times n$. Again below we exhibit a $3 \times 4 \times 2$ and a $3 \times 4 \times 3$ tour (which are equivalent to a $4 \times 3 \times 2$ and a $4 \times 3 \times 3$ tour)  with sites in the top left corners of the top and bottom layers.\\
\begin{center}
\includegraphics[scale=0.2]{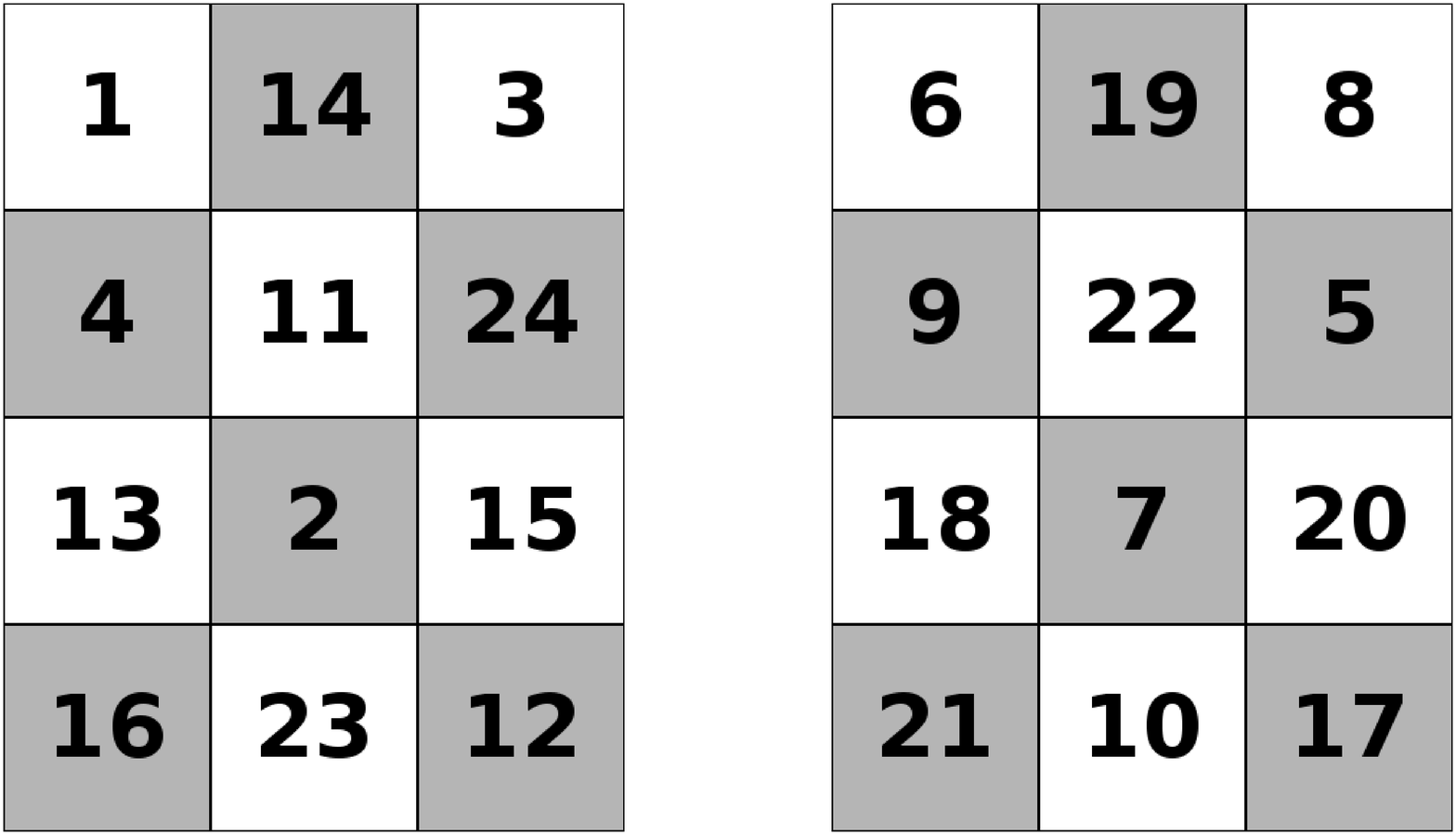}
\ \\
\ \\
\includegraphics[scale=0.2]{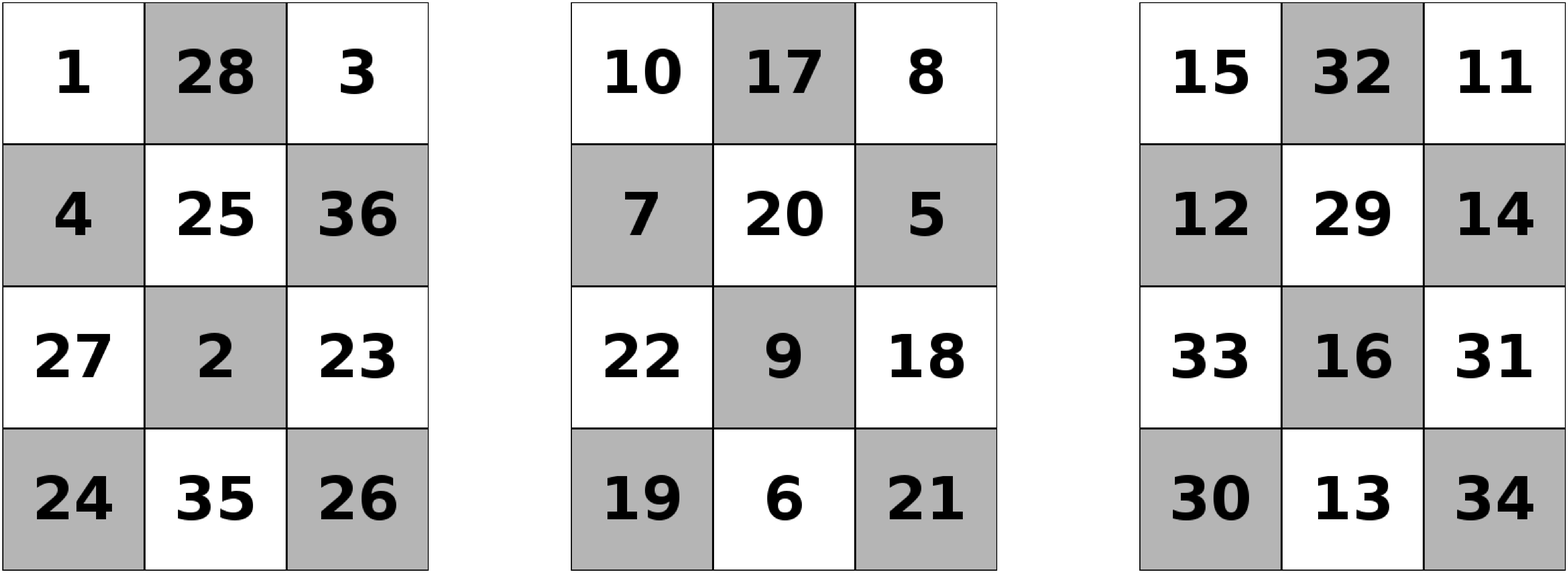}
\end{center}
\ \\
By the same method as above we can use these to construct bi-sited $4 \times 3 \times n$ tours for all $n \geq 2$.\\
\ \\
Finally we do $4 \times 2 \times n$. A $4 \times 2 \times 2$ tour does not exists, so we will need to have to use as our base cases a $4 \times 2 \times 3$, a $4 \times 2 \times 4$ and a $4 \times 2 \times 5$ tour.\\
\ \\
We can construct a $4 \times 6 \times 2$ tour by stacking two copies of the $3 \times 4 \times 2$ tour together, removing the $11-12$ line from the left copy and the $1-2$ line from the right copy and adding in the $11-11$ and $2-12$ lines. In this way we can construct bi-sited $4 \times n \times 2$ tours for all $n \equiv 0 $ mod$(3)$.\\
\ \\
Similarly we can add the $4 \times 3 \times 2$ tour to the left of the  $4 \times 4 \times 2$ tour we constructed by removing the $2-3$ line from the $4 \times 4 \times 2$ tour and the $1-2$ from the $4 \times 3 \times 2$ tour and adding in the $1-2$ and $2-3$ lines. This settles the case $n \equiv 1$ mod$(3)$.\\
\ \\
So we just need to construct a $5 \times 4 \times 2$ with two sites that includes the line $\big( (3,1,0)\,,\,(4,3,0) \big)$, which follows:\\
\begin{center}
\includegraphics[scale=0.2]{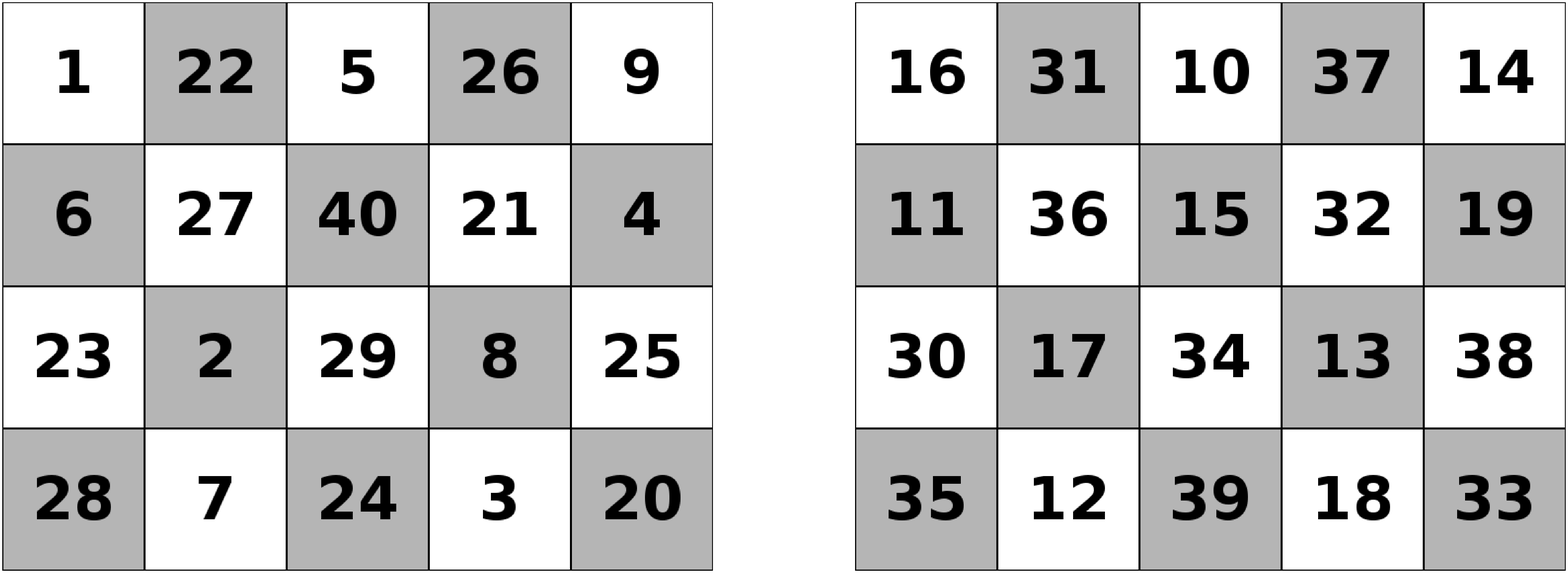}
\end{center}
\ \\
Now we look at the case $3 \times 2 \times n$ for $n \geq 4$. We can construct a $3 \times 2 \times 8$ tour by stacking together two copies of the $3 \times 2 \times 4$ tour, removing the line $15-16$ in the first copy and the line $8-9$ in the second copy and adding in the lines $15-8$ and $16-9$. By induction we can construct $3 \times 2 \times n$ tours for all $n \equiv 0$ mod$(4)$. \\
\ \\
So it will be sufficient to exhibit tours of size $3 \times 5 \times 2$, $3 \times 6 \times 2$ and $3 \times 7 \times 2$ with similar properties, which follow:\\
\begin{center}
\includegraphics[scale=0.2]{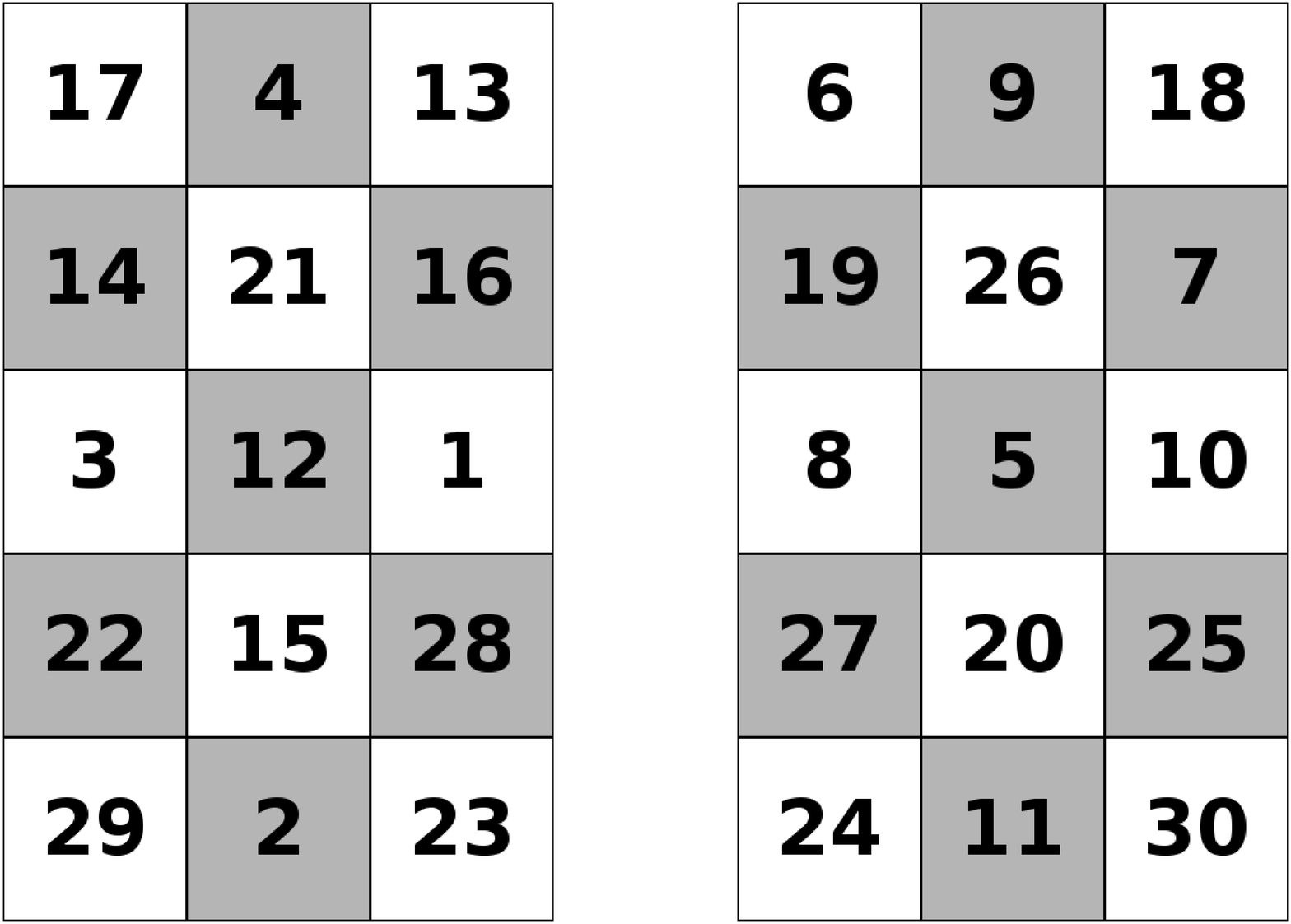}
\ \\
\ \\
\includegraphics[scale=0.2]{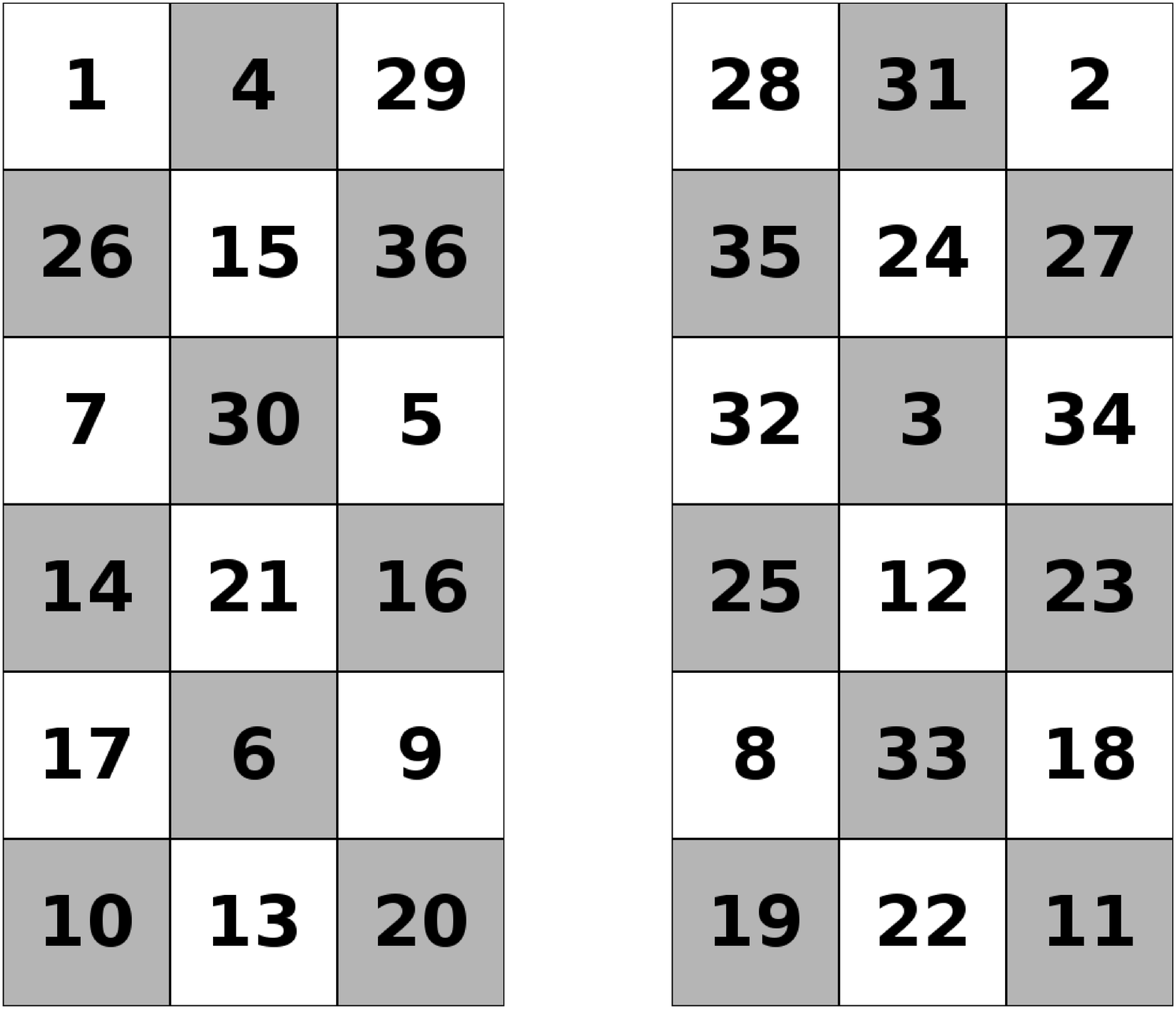}
\ \\
\ \\
\includegraphics[scale=0.2]{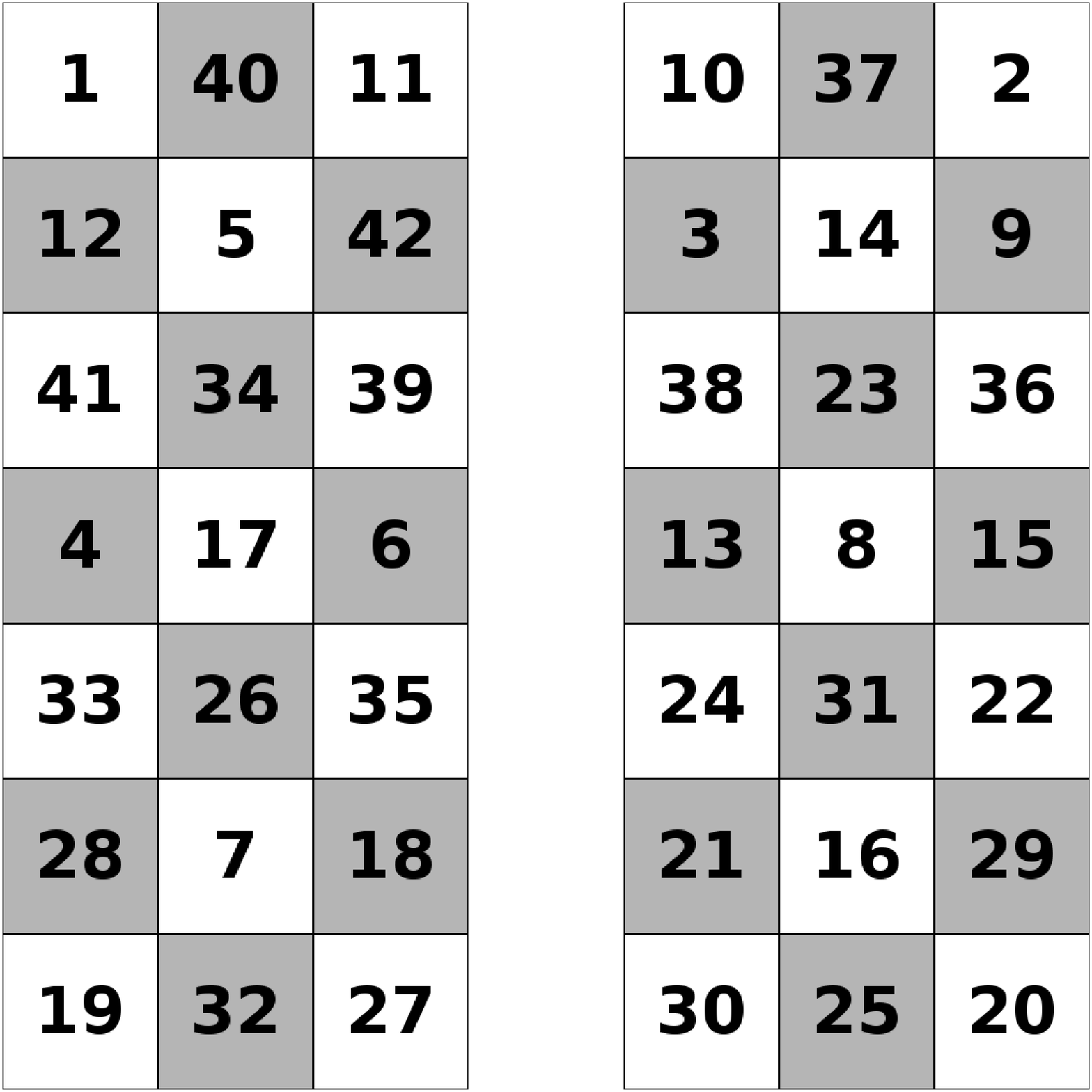}
\end{center}
\ \\
The remaining two cases are a $3 \times 3 \times 6$ and a $3 \times 3 \times 8$ tour. Firstly if we look back at the $4 \times 3 \times 3$ tour we can join two of these together to form a $8 \times 3 \times 3$ tour by deleting the $23-24$ line in the first copy and the $7-8$ line in the second copy and adding in the lines $7-24$ and $8-23$. The $3 \times 3 \times 6$ is below:\\
\begin{center}
\includegraphics[scale=0.2]{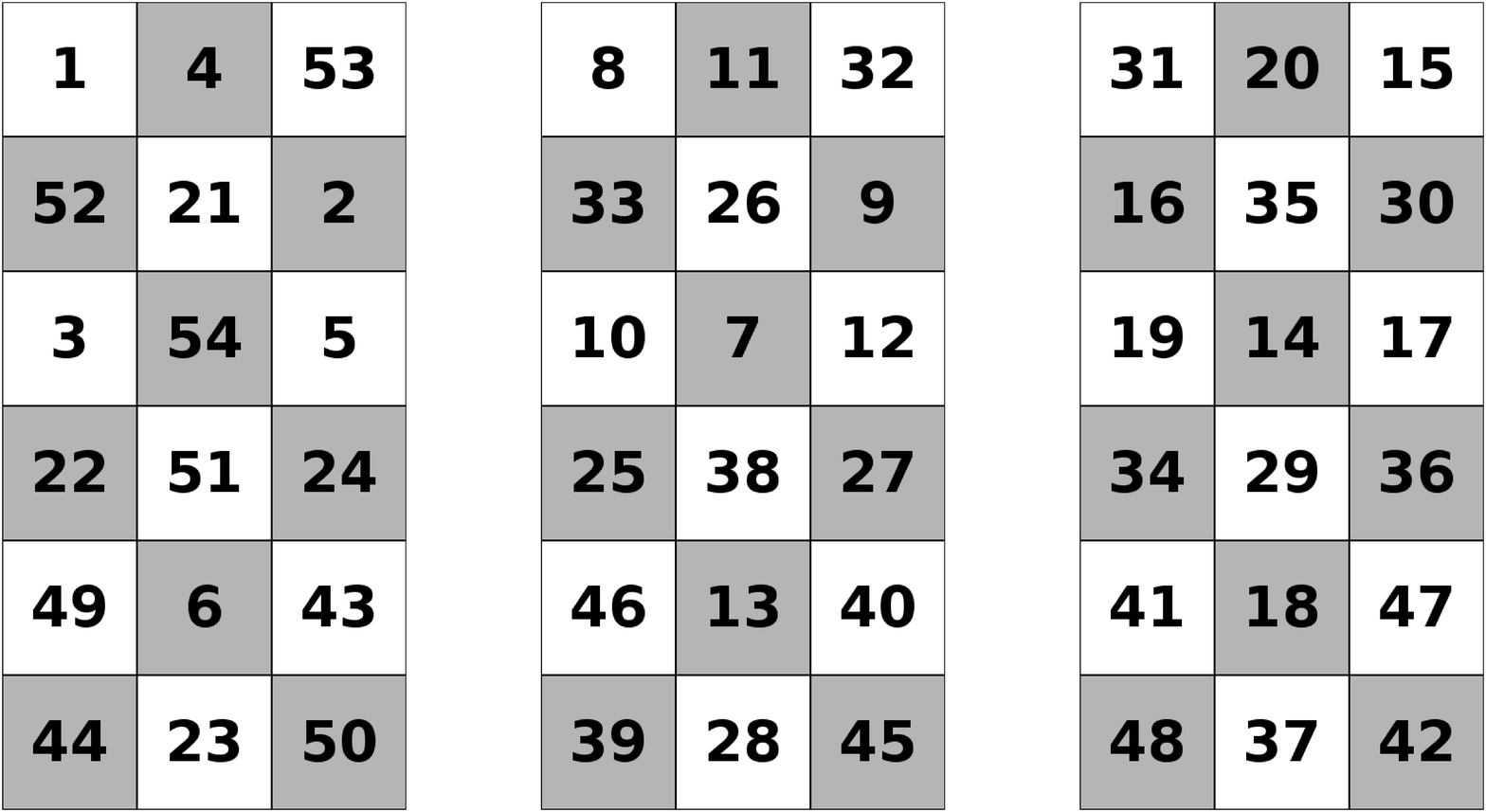}
\end{center}
where two sites are the lines $42-1$ and $11-12$ and the lines $2-3$ and $9-10$.\\
\bibliography{Knight}

\begin{thebibliography}{1}

\bibitem{D2007}
J.~DeMaio.
\newblock Which chessboards have a closed knight's tour within the cube?
\newblock {\em The Electronic Journal of Combinatorics}, 14, 2007.

\bibitem{DM2011}
J.~DeMaio and B.~Mathew.
\newblock Which chessboards have a closed knight's tour within the rectangular
  prism?
\newblock {\em The Electronic Journal of Combinatorics}, 18, 2011.

\bibitem{K1994}
D.~E. Knuth.
\newblock Leaper graphs.
\newblock {\em Mathematical Gazette}, 78:274--297, 1994.

\bibitem{S1991}
A.~J. Schwenk.
\newblock Which rectangular chessboards have a knight's tour?
\newblock {\em Mathematics Magazine}, 64(5):325--332, 1991.

\bibitem{S1971}
I.~Stewart.
\newblock Solid knight's tours.
\newblock {\em Journal of Recreational Mathematics}, 4(1), 1971.

\bibitem{W2004}
J.~J. Watkins.
\newblock {\em Across the Board: The Mathematics of Chessboard Problems}.
\newblock Princeton University Press, 2004.

\end{thebibliography}
\bibliographystyle{plain}
\end{document}